\documentclass[10pt,twoside]{article}

\usepackage{amsmath, amsthm, amssymb, mathrsfs, sansmath, setspace}

\newcommand{\R}		{\mathbb{R}}
\newcommand{\C}		{\mathbb{C}}

\newcommand{\Z}		{\mathbb{Z}}
\newcommand{\ds}	{\displaystyle}

\newtheorem{remark}{Remark}

\title{Backward Ricci Flow of Compact Locally Homogeneous Geometries on $4$-Manifolds}
\author{Thomas Bell}
\begin{document}

\bibliographystyle{plain}
\maketitle

\begin{abstract}
In this paper we study backward Ricci flow of locally homogeneous geometries of $4$-manifolds which admit compact quotients. We describe the long-term behavior of each class and show that many of the classes exhibit the same behavior near the singular time. In most cases, these manifolds converge to a sub-Riemannian geometry after suitable rescaling.
\end{abstract}

\section{Introduction}

Ricci flow on a manifold, $(M,g_0)$ is an evolution equation on the metric tensor given by

\begin{equation}\label{RF}
\frac{\partial g}{\partial t}=-2\text{Ric}_{g(t)}.
\end{equation}

This equation was first introduced by Richard Hamilton in \cite{Ha}, where he demonstrated that one can always expect short time existence as long as $(M,g_0)$ is a smooth manifold. Backward Ricci flow is described by 

\begin{equation}\label{BRF}
\frac{\partial g}{\partial t}=2\text{Ric}_{g(t)}.
\end{equation}

In general, we cannot expect short time existence of solutions to this equation. However, in the case of locally homogeneous manifolds, Ricci flow reduces to a system of Ordinary Differential Equations. As mentioned in \cite{CSC}, this eliminates all barriers to backward short-time existence of Ricci flow in these manifolds.\\

In \cite{IJ}, Isenberg and Jackson studied the behavior of solutions to Ricci flow along locally homogeneous $3$-manifolds. Later, in \cite{IJL}, Isenberg Jackson and  Lu studied Ricci flow along locally homogeneous $4$-manifolds which admit compact quotients. Subsequently, in \cite{CSC}, Cao and Saloff-Coste used calculations from \cite{CK} and \cite{IJ} to study backward Ricci flow on the homogeneous $3$-manifolds.\\

In this paper we use many of the calculations in \cite{IJL} to examine backward Ricci flow of compact locally homogeneous geometries on $4$-manifolds. The analysis in this paper is often very similar to that in \cite{IJL}, and we will use many of the same calculations, some of which are included for completeness.\\

The classes of locally homogeneous manifolds are described in \cite{IJL} and \cite{Ma}. For the geometries which are also Lie Groups we will choose a particular basis for the Lie Algebra, $\{X_1,X_2,X_3,X_4\}$, that satisfies certain bracket relations. More detail on these classes can be found in \cite{Ma}. Letting $\{\phi_i\}_{i=1}^4$ be the frame of $1$-forms dual to $\{X_i\}$ we can form a metric $g_0=A^{ij}\phi_i\otimes\phi_j$.\\

Our solutions will represent diagonalized Riemannian metrics. Thus they only exist as long as all unknowns remain positive and finite. We denote by $T_0$ the postiive time at which the various solutions to (\ref{BRF}) fail to exist.\\

We divide this paper into two main sections. First is the interesting section where we describe the Bianchi cases. The Lie Group structure has trivial Isotopy group, so the manifold is in itself a Lie Group. In the next section we give a quick description of each of the $4$-dimensional non-Bianchi cases. These are all metrics of constant sectional or holomorphic bi-sectional curvature or products of such metrics. The evolution of these metrics is well understood, but we include it in this paper for completeness. In the last section we summarize and compare the behaviors under backward Ricci flow of the various geometries under.

\section{The Bianchi cases}
The reduction of Ricci flow, equation (\ref{RF}), on the Bianchi classes of $4$-manifolds, to a system of ODE was done in \cite{IJL} using the following Ricci curvature formula for unimodular Lie groups from \cite[p.~184]{Be}. Recall that the elemenets $X,Y,Z,W$ come from the Lie Algebra $\mathfrak{g}$ of our Lie Group $G$ which represents our manifold.
\begin{equation}
Ric(W,W)=-\frac{1}{2}\sum_i\bigl|[W,Y_i]\bigr|^2-\frac{1}{2}\sum_i\bigl<\bigl[W,[W,Y_i]\bigr],Y_i\bigr>+\frac{1}{2}\sum_{i<j}\bigl<[Y_i,Y_j],W\bigr>^2.
\end{equation}
The only difference in calculating the system of ODE for backwards Ricci flow, equation (\ref{BRF}), is that the evolutions of the various metrics are negative of those found in \cite{IJL}. In this section we use these systems of ODE without further explanation.\\

A formula to calculate the sectional curvature on Lie Groups, with corresponding Lie Algebras, is found in \cite[p.~183]{Be}:

\begin{align}
\bigl<R(X,Y)X,Y\bigr>&=-\frac{3}{4}\bigl|[X,Y]\bigr|^2-\frac{1}{2}\bigl<\bigl[X,[X,Y]\bigr],Y\bigr>-\frac{1}{2}\bigl<\bigl[Y,[Y,X]\bigr],X\bigr>\notag\\
&\qquad\qquad +\bigl|U(X,Y)\bigr|^2-\bigl<U(X,X),U(Y,Y)\bigr>,\label{curvature}
\end{align}
where $U$ is defined by
\begin{equation}\label{U}
\bigl<U(X,Y),Z\bigr>=\frac{1}{2}\bigl<[Z,X],Y\bigr>+\frac{1}{2}\bigl<X,[Z,Y]\bigr>\text{ for all }Z\in \mathfrak{g}.
\end{equation}

The classification notation we use comes from \cite{IJL} and \cite{Ma}. See \cite{Ma} for more details on this classification.

\subsection{A1. Class $U1[(1,1,1)]$.}
Here we may choose a basis for the Lie Albegra $\{X_1,X_2,X_3,X_4\}$ such that the Lie bracket is of the form
\begin{align*}
[X_1,X_2]&=0 & [X_1,X_3]&=0 & [X_1,X_4]&=0\\
[X_2,X_3]&=0 & [X_2,X_4]&=0 & [X_3,X_4]&=0.
\end{align*}
Our Lie Group structure is $(M,G)=(\R^4,\R^4)$. The metric $g$ is flat, hence also Ricci flat and thus remains constant.
\subsection{A2. Class $U1[1,1,1]$.}
Here we may choose a basis for the Lie Albegra $\{X_1,X_2,X_3,X_4\}$ such that the Lie bracket is of the form
\begin{align*}
[X_1,X_2]&=0 & [X_1,X_3]&=0 & [X_1,X_4]&=X_1\\
[X_2,X_3]&=0 & [X_2,X_4]&=kX_2 & [X_3,X_4]&=-(k+1)X_3.
\end{align*}
When $k=0$ this corresponds to the geometry $(M,G)=\bigl(\tilde{Sol^3}\times\R,\tilde{Sol^3}\times\R\bigr)$.\\
When $k=1$ this corresponds to the geometry $(M,G)=\bigl(Sol_0^4,Sol_0^4\bigr)$.\\
If $k\neq 0,1$ and there is some number $\alpha>0$ with $\beta=k\alpha,~\gamma=-(k+1)\alpha$ such that $e^{\alpha},~e^{\beta}$ and $e^{\gamma}$ are roots of $\lambda^3-m\lambda^2+n\lambda-1=0$, then this corresponds to the geometry $(M,G)=\bigl(Sol^4_{m,n},Sol^4_{m,n}\bigr)$.\\

We diagonalize our initial metric $g_0$ by letting $Y_i=\Lambda_i^kX_k$ for constants $\Lambda_i^k$. Now letting $\{\theta_i\}$ be the frame of $1$-forms dual to $\{Y_i\}$ we may write the metric as
\begin{equation*}
g_0=\lambda_1\theta_1^2+\lambda_2\theta_2^2+\lambda_3\theta_3^2+\lambda_4\theta_4^2.
\end{equation*}

In \cite{IJL} we find outlines as to when exactly the Ricci tensor is also diagonal under these same coordinates. This is exactly when the metric will remain diagonal under Ricci flow and also backward Ricci flow. The property of a metric to remain diagonal under this flow is essential to our calculations, and we will only consider those families which satisfy this property. In each remaining subsection we will describe a flow which is of the form
\begin{equation*}
g(t)=A(t)\theta_1^2+B(t)\theta_2^2+C(t)\theta_3^2+D(t)\theta_4^2,
\end{equation*}
where
\begin{equation*}
A(0)=\lambda_1,~B(0)=\lambda_2,~C(0)=\lambda_3,\text{ and }D(0)=\lambda_4.
\end{equation*}

The systems of ODE governing the evolution of the quantities $A,~B,~C$ and $D$ were calculated in \cite{IJL}.\\

In the class $U1[(1,1,1)]$ we may diagonalize the metric by letting $Y_i=\Lambda_i^kX_k$ with
\begin{equation*}
\Lambda=\left[\begin{array}{cccc}
1&0&0&0\\a_1&1&0&0\\a_2&a_3&1&0\\a_4&a_5&a_6&1
\end{array}\right].
\end{equation*}
By Proposition 1 in \cite{IJL}, if $k=1$ then the metric $g(t)$ remains diagonal in the basis $Y_i$ if and only if $a_1=a_2=a_3=0$, and if $k\neq 1$ then $g(t)$ remains diagonal in the basis $Y_i$ if and only if $a_2=a_3=0$. In either of these cases we find that $\{Y_i\}$ satisfies the following bracket relations:
\begin{align*}
[Y_1,Y_2]&=0 & [Y_1,Y_3]&=0 & [Y_1,Y_4]&=Y_1\\
[Y_2,Y_3]&=0 & [Y_2,Y_4]&=kY_2 & [Y_3,Y_4]&=-(k+1)Y_3.
\end{align*}
Let $g(t)=A(t)\theta_1^2+B(t)\theta_2^2+C(t)\theta_3^2+D(t)\theta_4^2$ where $A(0)=\lambda_1,~B(0)=\lambda_2,~C(0)=\lambda_3$ and $D(0)=\lambda_4$. Then backward Ricci flow (\ref{BRF}) reduces to the following system of equations by \cite{IJL}:
\begin{equation}
\begin{aligned}
\frac{dA}{dt}&=\frac{dB}{dt}=\frac{dC}{dt}=0,\\
\frac{dD}{dt}&=-4\bigl(k^2+k+1\bigr).
\end{aligned}
\end{equation}
The solution is
\begin{equation}
\begin{aligned}
A(t)&=\lambda_1\\
B(t)&=\lambda_2\\
C(t)&=\lambda_3\\
D(t)&=\lambda_4-4\bigl(k^2+k+1\bigr)t.
\end{aligned}
\end{equation}

It is clear that
\begin{equation}\label{A2T0}
T_0=\frac{\lambda_4}{4}\bigl(k^2+k+1\bigr)^{-1},
\end{equation}
and that as $t\rightarrow T_0$ the volume normalized flow approaches the hyperplane, $\mathbb{R}^3$.\\

The sectional curvatures, also calculated in \cite{IJL}, are as follows:
\begin{align*}
K(Y_1,Y_2)&=-\frac{k}{D}& K(Y_1,Y_3)&=\frac{k+1}{D}& K(Y_1,Y_4)&=\frac{k(k+1)}{D}\\
K(Y_2,Y_3)&=-\frac{1}{D}& K(Y_2,Y_4)&=-\frac{k^2}{D}& K(Y_3,Y_4)&=-\frac{(k+1)^2}{D}.
\end{align*}
Thus, as $D$ approaches $0$ linearly in $t$ we see that the non-zero curvatures approach infinity with a singularity of the form $(T_0-t)^{-1}$.
\subsection{A3. Class $U1[\Z,\bar{\Z},1]$.}
Here we may choose a basis for the Lie Albegra $\{X_1,X_2,X_3,X_4\}$ such that the Lie bracket is of the form
\begin{align*}
[X_1,X_2]&=0 & [X_1,X_3]&=0 & [X_1,X_4]&=kX_1+X_2\\
[X_2,X_3]&=0 & [X_2,X_4]&=-X_1+kX_2c& [X_3,X_4]&=-2kX_3.
\end{align*}
This corresponds to the geometry $(M,G)=(\R^4,E(2)\times\R^2)$. 
Here we diagonalize the metric by letting $Y_i=\Lambda_i^kX_k$ with
\begin{equation*}
\Lambda=\left[\begin{array}{cccc}
1&a_2&a_3&0\\0&1&a_1&0\\0&0&1&0\\a_4&a_5&a_6&1
\end{array}\right].
\end{equation*}
By Proposition 2 in \cite{IJL}, the metric $g(t)$ remains diagonal in the basis $Y_i$ if and only if $a_1=a_2=a_3=0$. We find that $\{Y_i\}$ satisfies the following bracket relations:
\begin{align*}
[Y_1,Y_2]&=0 & [Y_1,Y_3]&=0 & [Y_1,Y_4]&=kY_1+Y_2\\
[Y_2,Y_3]&=0 & [Y_2,Y_4]&=-Y_1+kY_2 & [Y_3,Y_4]&=-2kY_3.
\end{align*}
Let $g(t)=A(t)\theta_1^2+B(t)\theta_2^2+C(t)\theta_3^2+D(t)\theta_4^2$ where $A(0)=\lambda_1,~B(0)=\lambda_2,~C(0)=\lambda_3$ and $D(0)=\lambda_4$. Then backward Ricci flow reduces to the following system of equations:
\begin{equation}\label{A3eqn}
\begin{aligned}
\frac{dA}{dt}&=\frac{A^2-B^2}{BD}\\
\frac{dB}{dt}&=\frac{B^2-A^2}{AD}\\
\frac{dC}{dt}&=0\\
\frac{dD}{dt}&=-\frac{(A-B)^2+12k^2AB}{AB}.
\end{aligned}
\end{equation}

Clearly, $C(t)=\lambda_3$. Also,
\begin{equation*}
\frac{d}{dt}\bigl(A-B\bigr)=\frac{(A+B)^2}{ABD}(A-B),
\end{equation*}
so the conditions $A=B,~A>B$ and $A<B$ are all preserved.\\

First assume $\lambda_1=\lambda_2$. Then $A=B$ for all time $t$, and our solution is

\begin{equation}
\begin{aligned}
A(t)&=\lambda_1\\
B(t)&=\lambda_1\\
C(t)&=\lambda_3\\
D(t)&=\lambda_4-12k^2t.
\end{aligned}
\end{equation}

Thus
\begin{equation}
T_0=\frac{\lambda_4}{12k^2},
\end{equation}
and it is clear that as $t\rightarrow T_0$, the volume normalized flow approaches the hyperplane, $\mathbb{R}^3$.\\

Now assume $\lambda_1\neq\lambda_2$. By the symmetry of (\ref{A3eqn}), we may assume that $\lambda_1>\lambda_2$. Then $A(t)\geq B(t)$ for all time $t$. Now, 
\begin{equation*}
\frac{d}{dt}\bigl(AB\bigr)=A\frac{B^2-A^2}{AD}+B\frac{A^2-B^2}{BD}=0,
\end{equation*}
so $AB=\lambda_1\lambda_2$ is constant. It is also clear that $T_0<\infty$, because $\ds\frac{dD}{dt}< -12k^2$ implies $T_0<\frac{\lambda_4}{12k^2}$. With $AB=\lambda_1\lambda_2$, we see that
\begin{equation*}
\frac{dA}{dt}=\frac{A^4-\lambda_1^2\lambda_2^2}{\lambda_1\lambda_2AD},~\frac{dD}{dt}=-\frac{\bigl(A-B\bigr)^2}{\lambda_1\lambda_2}-12k^2,
\end{equation*}
hence
\begin{equation}
\frac{1}{D}\frac{dD}{dA}=-\frac{\left[\left(A-\frac{\lambda_1\lambda_2}{A}\right)^2+12\lambda_1\lambda_2k^2\right]A}{A^4-\lambda_1^2\lambda_2^2}.
\end{equation}

Solving this equation we find
\begin{equation}
D=\Lambda\left(\frac{\sqrt{\lambda_1\lambda_2}A}{A^2+\lambda_1\lambda_2}\right)\left(\frac{A^2+\lambda_1\lambda_2}{A^2-\lambda_1\lambda_2}\right)^{3k^2},
\end{equation}
where
\begin{equation}\label{A3lambda}
~\Lambda=\lambda_4\left(\frac{\lambda_1+\lambda_2}{\sqrt{\lambda_1\lambda_2}}\right)\left(\frac{\lambda_1-\lambda_2}{\lambda_1+\lambda_2}\right)^{3k^2}.
\end{equation}

We see that $D\rightarrow 0$ as $t\rightarrow T_0$ if and only if $A\rightarrow\infty$ as $t\rightarrow T_0$. But since $AB$ is constant, we know that $B\rightarrow 0$ if and only if $A\rightarrow\infty$. Thus we know that $D\rightarrow 0,~B\rightarrow 0$ and $A\rightarrow\infty$ at the same time, $T_0< \frac{\lambda_4}{12k^2}<\infty$.\\

As $A\rightarrow\infty$, we see that
\begin{equation}\label{A3D}
D\approx \frac{\Lambda\sqrt{\lambda_1\lambda_2}}{A}=\frac{\Lambda}{\sqrt{\lambda_1\lambda_2}}B.
\end{equation}

Now, as $t\rightarrow T_0$, we may approximate
\begin{equation}\label{A3A1}
\frac{dA}{dt}=\frac{A^2-B^2}{BD}\sim\frac{A^4}{\Lambda(\lambda_1\lambda_2)^{3/2}}.
\end{equation}
Thus
\begin{equation*}
\frac{d}{dt}\left(\frac{1}{A^3}\right)\rightarrow -\frac{3}{\Lambda(\lambda_1\lambda_2)^{3/2}}
\end{equation*}
as $t\rightarrow T_0$. Since $A\rightarrow\infty$ as $t\rightarrow T_0$, we see that
\begin{equation*}
\frac{1}{A^3}=\frac{3}{\Lambda(\lambda_1\lambda_2)^{3/2}}(T_0-t)\bigl(1+o(T_0-t)\bigr).
\end{equation*}
Thus
\begin{align*}
A&=\left(\frac{3}{\Lambda(\lambda_1\lambda_2)^{3/2}}(T_0-t)\right)^{-1/3}\bigl(1+o(T_0-t)\bigr)^{-1/3}\\
&=\sqrt{\lambda_1\lambda_2}\left(\frac{3}{\Lambda}(T_0-t)\right)^{1/3}\bigl(1+o(T_0-t)\bigr)
\end{align*}

Now by (\ref{A3D}) we have the following behavior as $t\rightarrow T_0$:
\begin{equation}
\begin{aligned}
A&\approx\sqrt{\lambda_1\lambda_2}\left(\frac{3}{\Lambda}\bigl(T_0-t\bigr)\right)^{-1/3}\\
B&\approx\sqrt{\lambda_1\lambda_2}\left(\frac{3}{\Lambda}\bigl(T_0-t\bigr)\right)^{1/3}\\
C&=\lambda_3\\
D&\approx\Lambda\left(\frac{3}{\Lambda}\bigl(T_0-t\bigr)\right)^{1/3}.
\end{aligned}
\end{equation}
where $\Lambda$ is given by (\ref{A3lambda}).

The volume normalized solution will converge to the plane $\R^2$.\\

The sectional curvatures are as follows:
\begin{align*}
K(Y_1,Y_2)&=-\frac{\frac{A}{B}+\frac{B}{A}-2-4k^2}{4D}\hspace{.2in} & K(Y_1,Y_3)&=\frac{2k^2}{D}\\
K(Y_1,Y_4)&=\frac{\frac{A}{B}-3\frac{B}{A}+2-4k^2}{D}\hspace{.2in} & K(Y_2,Y_3)&=\frac{2k^2}{D}\\ K(Y_2,Y_4)&=\frac{-3\frac{A}{B}+\frac{B}{A}+2-4k^2}{D}\hspace{.2in} & K(Y_3,Y_4)&=-\frac{4k^2}{D}.
\end{align*}

Thus we see that the sectional curvatures perpendicular to $Y_3$ approach infinity at a rate of $(T_0-t)^{-1}$. If $k=0$ the sectional curvatures parallel to $Y_3$ remain $0$, while if $k\neq 0$ then these curvatures approach infinity at a rate of $(T_0-t)^{-1/3}$.
\subsection{A4. Class $U1[2,1],~\mu=0$.}
Here we may choose a basis for the Lie Albegra $\{X_1,X_2,X_3,X_4\}$ such that the Lie bracket is of the form
\begin{align*}
[X_1,X_2]&=0 & [X_1,X_3]&=0 & [X_1,X_4]&=X_2\\
[X_2,X_3]&=0 & [X_2,X_4]&=0 & [X_3,X_4]&=0.
\end{align*}
This corresponds to the geometry $(M,G)=(Nil^3\times\R,Nil^3\times\R)$. 
Here we diagonalize the metric by letting $Y_i=\Lambda_i^kX_k$ with
\begin{equation*}
\Lambda=\left[\begin{array}{cccc}
1&a_2&a_3&0\\0&1&0&0\\0&a_1&1&0\\a_4&a_5&a_6&1
\end{array}\right].
\end{equation*}
By Proposition 3 in \cite{IJL}, the metric $g(t)$ remains diagonal in the basis $Y_i$ for all $t$. We find that $\{Y_i\}$ satisfies the following bracket relations:
\begin{align*}
[Y_1,Y_2]&=0 & [Y_1,Y_3]&=0 & [Y_1,Y_4]&=Y_2\\
[Y_2,Y_3]&=0 & [Y_2,Y_4]&=-Y_1+kY_2 & [Y_3,Y_4]&=-2kY_3.
\end{align*}
Let $g(t)=A(t)\theta_1^2+B(t)\theta_2^2+C(t)\theta_3^2+D(t)\theta_4^2$ where $A(0)=\lambda_1,~B(0)=\lambda_2,~C(0)=\lambda_3$ and $D(0)=\lambda_4$. Then backward Ricci flow reduces to the following system of equations:
\begin{equation}\label{A4eqn}
\begin{aligned}
\frac{dA}{dt}&=-\frac{B}{D}\\
\frac{dB}{dt}&=\frac{B^2}{AD}\\
\frac{dC}{dt}&=0\\
\frac{dD}{dt}&=-\frac{B}{A}.
\end{aligned}
\end{equation}

We calculate
\begin{equation*}
\frac{d}{dt}(AB)=\frac{d}{dt}\left(\frac{A}{D}\right)=0.
\end{equation*}

Thus we have 
\begin{equation}\label{A4BandD}
B=\frac{\lambda_1\lambda_2}{A},\text{ and }D=\frac{\lambda_4}{\lambda_1}A,
\end{equation}
so
\begin{equation}
A^2\frac{dA}{dt}=-AB\frac{A}{D}=-\frac{\lambda_1^2\lambda_2}{\lambda_4},
\end{equation}
hence
\begin{equation}\label{A4A}
A^3=\lambda_1^3-\frac{3\lambda_1^2\lambda_2}{\lambda_4}t.
\end{equation}

Thus by (\ref{A4BandD}) and (\ref{A4A}), we have the following solution to (\ref{A4eqn}):
\begin{equation}
\begin{aligned}
A&=\lambda_1\left(1-\frac{3\lambda_2}{\lambda_1\lambda_4}t\right)^{1/3}\\
B&=\lambda_2\left(1-\frac{3\lambda_2}{\lambda_1\lambda_4}t\right)^{-1/3}\\
C&=\lambda_3\\
D&=\lambda_4\left(1-\frac{3\lambda_2}{\lambda_1\lambda_4}t\right)^{1/3}.
\end{aligned}
\end{equation}

We see that
\begin{equation}\label{A4T0}
T_0=\frac{\lambda_1\lambda_4}{3\lambda_2},
\end{equation}
and that the volume normalized solution will converge to the plane $\R^2$.\\

The sectional curvatures are as follows:
\begin{align*}
K(Y_1,Y_2)&=\frac{B}{4AD}& K(Y_1,Y_3)&=0 & K(Y_1,Y_4)&=-\frac{3B}{4AD}\\
K(Y_2,Y_3)&=0& K(Y_2,Y_4)&=\frac{B}{4AD} & K(Y_3,Y_4)&=0.
\end{align*}
Thus all non-zero curvatures approach infinity near $t=T_0$ at a rate of $(T_0-t)^{-1}$.

\subsection{A5. Class $U1[2,1],~\mu=1$}
Here we may choose a basis for the Lie Albegra $\{X_1,X_2,X_3,X_4\}$ such that the Lie bracket is of the form
\begin{align*}
[X_1,X_2]&=0 & [X_1,X_3]&=0 & [X_1,X_4]&=-\frac{1}{2}X_1+X_2\\
[X_2,X_3]&=0 & [X_2,X_4]&=-\frac{1}{2}X_2 & [X_3,X_4]&=X_3.
\end{align*}

This does not correspond to any of the compact homogeneous geometries. 
Here we diagonalize the metric by letting $Y_i=\Lambda_i^kX_k$ with
\begin{equation*}
\Lambda=\left[\begin{array}{cccc}
1&a_2&a_3&0\\0&1&a_1&0\\0&0&1&0\\a_4&a_5&a_6&1
\end{array}\right].
\end{equation*}
By Proposition 4 in \cite{IJL}, the metric $g(t)$ remains diagonal in the basis $Y_i$ if and only if $a_1=a_3=0$. We find that $\{Y_i\}$ satisfies the following bracket relations:
\begin{align*}
[Y_1,Y_2]&=0 & [Y_1,Y_3]&=0 & [Y_1,Y_4]&=-\frac{1}{2}Y_1+Y_2\\
[Y_2,Y_3]&=0 & [Y_2,Y_4]&=-\frac{1}{2}Y_2 & [Y_3,Y_4]&=Y_3.
\end{align*}

Let $g(t)=A(t)\theta_1^2+B(t)\theta_2^2+C(t)\theta_3^2+D(t)\theta_4^2$ where $A(0)=\lambda_1,~B(0)=\lambda_2,~C(0)=\lambda_3$ and $D(0)=\lambda_4$. Then backward Ricci flow reduces to the following system of equations:
\begin{equation}\label{A5eqn}
\begin{aligned}
\frac{dA}{dt}&=-\frac{B}{D}\\
\frac{dB}{dt}&=\frac{B^2}{AD}\\
\frac{dC}{dt}&=0\\
\frac{dD}{dt}&=-3-\frac{B}{A}.
\end{aligned}
\end{equation}

Since $\ds\frac{dD}{dt}<-3$ we see that there is a maximal time $\ds T_0<\frac{\lambda_4}{3}$. Also,
\begin{equation}
\frac{d}{dt}(AB)=0,
\end{equation}
so
\begin{equation}\label{A5AB}
AB=\lambda_1\lambda_2.
\end{equation}
Thus (\ref{A5eqn}) reduces to
\begin{equation}\label{A5eqn2}
\begin{aligned}
\frac{dA}{dt}&=-\frac{\lambda_1\lambda_2}{AD}\\
B&=\frac{\lambda_2}{\lambda_1}A\\
C&=\lambda_3\\
\frac{dD}{dt}&=\frac{-3A^2-\lambda_1\lambda_2}{A^2}.\\
\end{aligned}
\end{equation}

Now we calculate
\begin{align}
\frac{1}{D}\cdot\frac{dD}{dA}&=\frac{3}{\lambda_1\lambda_2}A+\frac{1}{A}.\notag
\intertext{Solving gives us}
D&=\Lambda\cdot Ae^{\left(\frac{3A^2}{2\lambda_1\lambda_2}\right)},\label{A5D}
\intertext{where}
\Lambda&=\frac{\lambda_4}{\lambda_1}e^{\left(-\frac{3\lambda_1}{2\lambda_2}\right)}.\label{A5lambda}
\end{align}

By (\ref{A5AB}) and (\ref{A5D}) we see that $A\rightarrow 0,~B\rightarrow\infty$ and $D\rightarrow 0$ as $t\rightarrow T_0$.\\

To describe the behavior near $t=T_0$ we observe that as $A$ approaches $0$, (\ref{A5eqn}), (\ref{A5AB}) and (\ref{A5D}) tell us
\begin{equation*}
A^2\frac{dA}{dt}\rightarrow -\frac{\lambda_1\lambda_2}{\Lambda},
\end{equation*}
hence
\begin{equation}
A^3=\frac{3\lambda_1\lambda_2}{\Lambda}(T_0-t)\bigl(1+o(T_0-t)\bigr).
\end{equation}

Thus, by  (\ref{A5eqn}), (\ref{A5AB}) and (\ref{A5D}), we have the following solutions to (\ref{A5eqn}):
\begin{equation}
\begin{aligned}
A&\approx\left(\frac{3\lambda_1\lambda_2}{\Lambda}\bigl(T_0-t\bigr)\right)^{1/3}\\
B&\approx\lambda_1\lambda_2\left(\frac{3\lambda_1\lambda_2}{\Lambda}\bigl(T_0-t\bigr)\right)^{-1/3}\\
C&=\lambda_3\\
D&\approx \Lambda\left(\frac{3\lambda_1\lambda_2}{\Lambda}\bigl(T_0-t\bigr)\right)^{1/3}.
\end{aligned}
\end{equation}
where $\Lambda$ is given in (\ref{A5lambda}). Again we notice that as $t$ approaches $T_0$ the renormalized flow approaches the plane $\R^2$.\\

The sectional curvatures are as follows:
\begin{align*}
K(Y_1,Y_2)&=\frac{-1+\frac{B}{A}}{4D} & K(Y_1,Y_3)&=\frac{1}{2D} & K(Y_1,Y_4)&=-\frac{1+3\frac{B}{A}}{4D}\\
K(Y_2,Y_3)&=\frac{1}{2D} & K(Y_2,Y_4)&=\frac{-1+\frac{B}{A}}{4D} & K(Y_3,Y_4)&=-\frac{1}{D}.
\end{align*}

Thus curvatures perpendicular to $Y_3$ will have a singularity $t=T_0$ of the form $(T_0-t)^{-1}$, while those parallel to $Y_3$ will have a singularity of the form $(T_0-t)^{-1/3}$.\\

\subsection{A6. Class U1[3].}
Here we may choose a basis for the Lie Albegra $\{X_1,X_2,X_3,X_4\}$ such that the Lie bracket is of the form
\begin{align*}
[X_1,X_2]&=0 & [X_1,X_3]&=0 & [X_1,X_4]&=X_2\\
[X_2,X_3]&=0 & [X_2,X_4]&=X_3 & [X_3,X_4]&=0.
\end{align*}
This corresponds to the geometry $(M,G)=(Nil^4,Nil^4)$.\\

Here we diagonalize the metric by letting $Y_i=\Lambda_i^kX_k$ with
\begin{equation*}
\Lambda=\left[\begin{array}{cccc}
1&a_2&a_3&0\\0&1&0&0\\0&a_1&1&0\\a_4&a_5&a_6&1
\end{array}\right].
\end{equation*}
By Proposition 5 in \cite{IJL}, the metric $g(t)$ remains diagonal in the basis $Y_i$ if and only if $a_1=a_2$. We find in this case that $\{Y_i\}$ satisfies the following bracket relations:
\begin{align*}
[Y_1,Y_2]&=0 & [Y_1,Y_3]&=0 & [Y_1,Y_4]&=Y_2\\
[Y_2,Y_3]&=0 & [Y_2,Y_4]&=Y_3 & [Y_3,Y_4]&=0.
\end{align*}
Let $g(t)=A(t)\theta_1^2+B(t)\theta_2^2+C(t)\theta_3^2+D(t)\theta_4^2$ where $A(0)=\lambda_1,~B(0)=\lambda_2,~C(0)=\lambda_3$ and $D(0)=\lambda_4$. Then backward Ricci flow reduces to the following system of equations:
\begin{equation}\label{A6eqn}
\begin{aligned}
\frac{dA}{dt}&=-\frac{B}{D}\\
\frac{dB}{dt}&=\frac{B^2-AC}{AD}\\
\frac{dC}{dt}&=\frac{C^2}{BD}\\
\frac{dD}{dt}&=-\frac{B}{A}-\frac{C}{B}.
\end{aligned}
\end{equation}

Similar Calculations as those in \cite{IJL} give us the following solution to (\ref{A6eqn}):
\begin{align}
&\begin{aligned}
A&=\lambda_1\left(1-\frac{3\lambda_2}{\lambda_1\lambda_4}t\right)^{1/3}\\
B&=\lambda_2\left(1-\frac{3\lambda_2}{\lambda_1\lambda_4}t\right)^{-1/3}\left(1-\frac{3\lambda_3}{\lambda_2\lambda_4}t\right)^{1/3}\\
C&=\lambda_3\left(1-\frac{3\lambda_3}{\lambda_2\lambda_4}t\right)^{-1/3}\\
D&=\lambda_4\left(1-\frac{3\lambda_2}{\lambda_1\lambda_4}t\right)^{1/3}\left(1-\frac{3\lambda_3}{\lambda_2\lambda_4}t\right)^{1/3}.
\end{aligned}
\intertext{ where}
&T_0=\min\left\{\frac{\lambda_1\lambda_4}{3\lambda_2},\frac{\lambda_2\lambda_4}{3\lambda_3}\right\}.
\end{align}

The sectional curvatures are as follows:
\begin{align*}
K(Y_1,Y_2)&=\frac{B}{4AD} & K(Y_1,Y_3)&=0 & K(Y_1,Y_4) &=\frac{-3B}{4AD}\\
K(Y_2,Y_3)&=\frac{C}{4BD} & K(Y_2,Y_4)&=\frac{\frac{B}{A}-3\frac{C}{B}}{4D} & K(Y_3,Y_4)&=\frac{C}{4BD}.
\end{align*}

If $\lambda_2^2<\lambda_1\lambda_3$, then $\ds T_0=\frac{\lambda_2\lambda_4}{3\lambda_3}$, and as $t\rightarrow T_0$ we have
\begin{equation}
\begin{aligned}
A&\approx k_1\\
B&\approx k_2(T_0-t)^{1/3}\\
C&=k_3(T_0-t)^{-1/3}\\
D& \approx k_4(T_0-t)^{1/3}.\\
\end{aligned}
\end{equation}

Also, $K(Y_1,Y_2)$ approaches a positive constant, $K(Y_1,Y_4)$ approaches a negative constant, and all curvatures perpendicular to $Y_1$ will have singularities of the form $(T_0-t)^{-1}$.\\

If $\lambda_2^2>\lambda_1\lambda_3$, then $\ds T_0=\frac{\lambda_1\lambda_4}{3\lambda_2}$, and as $t\rightarrow T_0$ we have
\begin{equation}
\begin{aligned}
A&= k_1(T_0-t)^{1/3}\\
B&\approx k_2(T_0-t)^{-1/3}\\
C&\approx k_3\\
D& \approx k_4(T_0-t)^{1/3}.\\
\end{aligned}
\end{equation}

Here $K(Y_2,Y_3)$ and $K(Y_3,Y_4)$ approach positive constants while each curvature perpendicular to $Y_3$ will have a singularity of the form $(T_0-t)^{-1}$.\\

If $\lambda_2^2=\lambda_1\lambda_3$, then $\ds T_0=\frac{\lambda_2\lambda_4}{3\lambda_3}=\frac{\lambda_1\lambda_4}{3\lambda_2}$, and as $t\rightarrow T_0$ we have
\begin{equation}
\begin{aligned}
A&= k_1(T_0-t)^{1/3}\\
B&=\lambda_2\\
C&=k_3(T_0-t)^{-1/3}\\
D&=k_4(T_0-t)^{2/3}.\\
\end{aligned}
\end{equation}

Here all non-zero curvatures will have a singularity of the form $(T_0-t)^{-1}$.\\

In all three cases the volume normalized solution approaches the plane $\R^2$.\\

\subsection{A7. Class $U3I0$.}
Here we may choose a basis for the Lie Albegra $\{X_1,X_2,X_3,X_4\}$ such that the Lie bracket is of the form
\begin{align*}
[X_1,X_2]&=-X_3 & [X_1,X_3]&=-X_2 & [X_1,X_4]&=0\\
[X_2,X_3]&=X_4 & [X_2,X_4]&=0 & [X_3,X_4]&=0.
\end{align*}
This corresponds to the geometry $(M,G)=(Sol^4,Sol^4)$. 
Here we diagonalize the metric by letting $Y_i=\Lambda_i^kX_k$ with
\begin{equation*}
\Lambda=\left[\begin{array}{cccc}
1&a_4&a_5&a_6\\0&1&a_2&a_3\\0&0&1&a_1\\0&0&0&1
\end{array}\right].
\end{equation*}
Let $\alpha=a_2,~\beta=a_1a_2-_3-a_4$ and $\gamma=a_1-a_1a_2^2+a_2a_3+a_2a_4-a_5$. Proposition 6 in \cite{IJL}, the metric $g(t)$ remains diagonal in the basis $Y_i$ if and only if one of the following hold:
\begin{center}
\begin{tabular}{ll}
(i) & $\alpha=\beta=\gamma=0$\\
(ii) & $\beta=\gamma=0$ and $\lambda_2=(1-\alpha^2)\lambda_3$
\end{tabular}
\end{center}
We analyze these cases separately.

\subsubsection{A7(i).}
Here we consider the case where $\alpha=\beta=\gamma=0$. We find that $\{Y_i\}$ satisfies the following bracket relations:
\begin{align*}
[Y_1,Y_2]&=-Y_3 & [Y_1,Y_3]&=-Y_2 & [Y_1,Y_4]&=0\\
[Y_2,Y_3]&=Y_4 & [Y_2,Y_4]&=0 & [Y_3,Y_4]&=0.
\end{align*}
Let $g(t)=A(t)\theta_1^2+B(t)\theta_2^2+C(t)\theta_3^2+D(t)\theta_4^2$ where $A(0)=\lambda_1,~B(0)=\lambda_2,~C(0)=\lambda_3$ and $D(0)=\lambda_4$. Then backward Ricci flow reduces to the following system of equations:
\begin{equation}\label{A7ieqn}
\begin{aligned}
\frac{dA}{dt}&=-\frac{B}{C}-\frac{C}{B}-2\\
\frac{dB}{dt}&=-\frac{C}{A}-\frac{D}{C}+\frac{B^2}{AC}\\
\frac{dC}{dt}&=-\frac{B}{A}-\frac{D}{B}+\frac{C^2}{AB}\\
\frac{dD}{dt}&=\frac{D^2}{BC}.
\end{aligned}
\end{equation}

By the symmetry of $B$ and $C$ in (\ref{A7ieqn}), we may assume that $\lambda_2\geq\lambda_3$. We calculate
\begin{equation*}
\frac{d}{dt}\bigl(BCD^2\bigr)=0,
\end{equation*}
so
\begin{equation}
BCD^2=\lambda_2\lambda_3\lambda_4^2.\label{A7iBCD}
\end{equation}

Now we observe
\begin{equation*}
\frac{1}{D^4}\cdot\frac{dD}{dt}=\frac{1}{\lambda_2\lambda_3\lambda_4^2},
\end{equation*}
and
\begin{equation}
D=\lambda_4\left(1-\frac{3\lambda_4}{\lambda_2\lambda_3}t\right)^{-1/3}.\label{A7iD}
\end{equation}

We may also calculate
\begin{equation*}
\frac{d}{dt}\bigl(AD(B-C)\bigr)=0
\end{equation*}
so 
\begin{equation}
AD(B-C)=\lambda_1\lambda_4(\lambda_2-\lambda_3).\label{A7iADBC}
\end{equation}
Now by (\ref{A7iBCD}) and (\ref{A7iADBC}) we have
\begin{equation}{\label{A7ik}}
\frac{(B-C)^2}{BC}A^2=\frac{\bigl(AD(B-C)\bigr)^2}{BCD^2}=\frac{\lambda_1^2(\lambda_2-\lambda_3)^2}{\lambda_2\lambda_3}=4k^2.
\end{equation}

Now
\begin{equation}
\frac{dA}{dt}=-\frac{(B+C)^2}{BC}=-\frac{4\bigl(k^2+A^2\bigr)}{A^2},
\end{equation}
hence
\begin{equation}A-k\tan^{-1}\left(\frac{A}{k}\right)=-4t+\lambda_1-k\tan^{-1}\left(\frac{\lambda_1}{k}\right),\label{A7iA}
\end{equation}
where $k$ is given in (\ref{A7ik}). Now by (\ref{A7iD}) and (\ref{A7iA}) we have
\begin{align}
T_0&=\min\left\{\frac{\lambda_2\lambda_3}{3\lambda_4},~\frac{\lambda_1}{4}-\frac{\lambda_1(\lambda_2-\lambda_3)}{2\sqrt{\lambda_2\lambda_3}}\tan^{-1}\left(\frac{2\sqrt{\lambda_2\lambda_3}}{\lambda_2-\lambda_3}\right)\right\}\notag\\
&=\min\bigl\{T_1,T_2\bigr\}.
\end{align}

We calculate
\begin{equation*}
\lim_{A\rightarrow 0}\left(\frac{A-k\tan^{-1}\bigl(\frac{A}{k}\bigr)}{A^3}\right)=\frac{1}{3k^2}.
\end{equation*}

Thus by (\ref{A7iA}), near $A=0$ we have
\begin{equation}
A\approx\bigl(12k^2(T_2-t)\bigr)^{1/3}.
\end{equation}

Using (\ref{A7iBCD}) and (\ref{A7iADBC}) we can calculate
\begin{align}
B&=\frac{\lambda_4}{2AD}\biggl(\lambda_1(\lambda_2-\lambda_3)+\sqrt{\lambda_1^2(\lambda_2-\lambda_3)^2+4A^2\lambda_2\lambda_3}\biggr),\label{A7iB}\\
C&=\frac{\lambda_4}{2AD}\biggl(\lambda_1(\lambda_3-\lambda_2)+\sqrt{\lambda_1^2(\lambda_2-\lambda_3)^2+4A^2\lambda_2\lambda_3}\biggr).\label{A7iC}
\end{align}

The sectional curvatures are as follows:
\begin{align*}
K(Y_1,Y_2)&=\frac{\frac{B}{C}-3\frac{C}{B}-2}{4A} & K(Y_1,Y_3)&=\frac{\frac{C}{B}-3\frac{B}{C}-2}{4A}\\
K(Y_1,Y_4)&=0 & K(Y_2,Y_3)&=\frac{B^2+C^2+2BC-3AD}{4ABC}\\
K(Y_2,Y_4)&=\frac{D}{4BC} & K(Y_3,Y_4)&=\frac{D}{4BC}.
\end{align*}

If $T_0=T_1<T_2$, then we have the following behavior as $t\rightarrow T_0$:
\begin{equation}
\begin{aligned}
A&\approx k_1\\
B&\approx k_2(T_0-t)^{1/3}\\
C&\approx k_3(T_0-t)^{1/3}\\
D&=k_4(T_0-t)^{-1/3}.
\end{aligned}
\end{equation}

In this case the curvatures parallel to $Y_1$ approach constants while those perpendicular to $Y_1$ will have a singularity of the form $(T_0-t)^{-1}$. The volume normalized solution approaches $\R^2$.\\

If $T_0=T_2<T_1$ and $\lambda_2=\lambda_3$, then as $t\rightarrow T_0$,
\begin{equation}
\begin{aligned}
A&\approx k_1(T_0-t)^{1/3}\\
B=C&\approx k_2\\
D&\approx k_4.
\end{aligned}
\end{equation}

Here the curvatures perpendicular to $Y_4$ approach infinity at a rate of $(T_0-t)^{-1/3}$ while those curvatures parallel to $Y_4$ approach constants. The normalized solution approaches $\R^3$.\\

If $T_0=T_2<T_1$ and $\lambda_2>\lambda_3$, then as $t\rightarrow T_0$,
\begin{equation}
\begin{aligned}
A&\approx k_1(T_0-t)^{1/3}\\
B&\approx k_2(T_0-t)^{-1/3}\\
C&\approx k_3(T_0-t)^{1/3}\\
D&\approx k_4.
\end{aligned}
\end{equation}

In this case the curvatures perpendicular to $Y_4$ approach infinity at a rate of $(T_0-t)^{-1}$ while those curvatures parallel to $Y_4$ approach constants. The normalized solution approaches $\R^3$.\\

If $T_0=T_1=T_2$ and $\lambda_2=\lambda_3$, then as $t\rightarrow T_0, AD\rightarrow$ constant, so
\begin{equation}
\begin{aligned}
A&\approx k_1(T_0-t)^{1/3}\\
B=C&\approx k_2(T_0-t)^{1/3}\\
D&= k_4(T_0-t)^{-1/3}.
\end{aligned}
\end{equation}
Here, $K(Y_1,Y_2)$ and $K(Y_1,Y_3)$ will have singularities at $t=T_0$ of the form $(T_0-t)^{-1/3}$. All other non-zero sectional curvatures will have singularities of the form $(T_0-t)^{-1}$. The normalized solution approaches $\R$.\\

If $T_0=T_1=T_2$ and $\lambda_2>\lambda_3$, then as $t\rightarrow T_0, AD\rightarrow$ constant, so
\begin{equation}
\begin{aligned}
A&\approx k_1(T_0-t)^{1/3}\\
B&\approx k_2\\
C&\approx k_3(T_0-t)^{2/3}\\
D&= k_4(T_0-t)^{-1/3}.
\end{aligned}
\end{equation}
Here, all non-zero curvatures will have singularities of the form $(T_0-t)^{-1}$. The normalized solution approaches $\R^2$.

\subsubsection{A7(ii)} 
Here we consider the case where $\beta=\gamma=0$ and $\lambda_2=(1-\alpha^2)\lambda_3$. We find that $\{Y_i\}$ satisfies the following bracket relations:
\begin{align*}
[Y_1,Y_2]&=-\alpha Y_2 & [Y_1,Y_3]&=\alpha Y_3-Y_2 & [Y_1,Y_4]&=0\\
[Y_2,Y_3]&=Y_4 & [Y_2,Y_4]&=0 & [Y_3,Y_4]&=0.
\end{align*}
Let $g(t)=A(t)\theta_1^2+B(t)\theta_2^2+C(t)\theta_3^2+D(t)\theta_4^2$ where $A(0)=\lambda_1,~B(0)=\lambda_2,~C(0)=\lambda_3$ and $D(0)=\lambda_4$. Then backward Ricci flow reduces to the following system of equations:
\begin{equation}\label{A7iieqn}
\begin{aligned}
\frac{dA}{dt}&=-\frac{B^2+2(1+\alpha^2)BC+(1-\alpha^2)^2C^2}{BC}\\
\frac{dB}{dt}&=\frac{-AD+B^2-(1-\alpha^2)^2C^2}{AC}\\
\frac{dC}{dt}&=\frac{-AD-B^2+(1-\alpha^2)^2C^2}{AB}\\
\frac{dD}{dt}&=\frac{D^2}{BC}\\
\lambda_2&=(1-\alpha^2)\lambda_3.
\end{aligned}
\end{equation}

We observe
\begin{equation*}
\frac{d}{dt}\bigl(B-(1-\alpha^2)C\bigr)=0,
\end{equation*}
so $B=(1-\alpha^2)C$ is preserved under (\ref{A7iieqn}), which then reduces to
\begin{equation}
\begin{aligned}
\frac{dA}{dt}&=-4\\
\frac{dB}{dt}&=-(1-\alpha^2)\frac{D}{B}\\
B&=(1-\alpha^2)C\\
\frac{dD}{dt}&=(1-\alpha^2)\frac{D^2}{B^2}.
\end{aligned}
\end{equation}
Now
\begin{equation*}
\frac{d}{dt}\bigl(BD\bigr)=0,
\end{equation*}
so
\begin{equation}\label{A7iiBD}
BD=\lambda_2\lambda_4.
\end{equation}
Thus
\begin{equation*}
\frac{dB}{dt}=-(1-\alpha^2)\frac{\lambda_2\lambda_4}{B^2},
\end{equation*}
so
\begin{equation*}
\frac{1}{3}B^3=\frac{1}{3}\lambda_2^3-(1-\alpha^2)\lambda_2\lambda_4t,
\end{equation*}
and we have the following solution to (\ref{A7iieqn}):
\begin{equation}\label{A7iisol}
\begin{aligned}
A&=\lambda_1-4t\\
B&=\bigl(\lambda_2^3-3(1-\alpha^2)\lambda_2\lambda_4t\bigr)^{1/3}\\
C&=\frac{1}{1-\alpha^2}\bigl(\lambda_2^3-3(1-\alpha^2)\lambda_2\lambda_4t\bigr)^{1/3}\\
D&=\lambda_2\lambda_4\bigl(\lambda_2^3-3(1-\alpha^2)\lambda_2\lambda_4t\bigr)^{-1/3},
\end{aligned}
\end{equation}
where
\begin{equation}\label{A7iiT0}
T_0=\min\left\{\frac{\lambda_1}{4},~\frac{\lambda_2^2}{3(1-\alpha^2)\lambda_4}\right\}=\min\bigl\{T_1,T_2\bigr\}.
\end{equation}

The sectional curvatures are as follows:
\begin{align*}
K(Y_1,Y_2)&=-\frac{1}{A} & K(Y_1,Y_3)&=-\frac{1}{A} & K(Y_1,Y_4)&=0\\
K(Y_2,Y_3)&=\frac{4BC-3AD}{4ABC} & K(Y_2,Y_4)&=\frac{D}{4BC} & K(Y_3,Y_4)&=\frac{D}{4BC}.
\end{align*}

If $T_0=T_1<T_2$, then near $t=T_0$ we have
\begin{equation}
\begin{aligned}
A&=4(T_0-t)\\
B&\approx k_2\\
C&\approx k_3\\
D&\approx k_4.
\end{aligned}
\end{equation}

In this case all sectional curvatures perpendicular to $Y_4$ will have a singularity at $T_0$ of the form $(T_0-t)^{-1}$ while those parallel to $Y_4$ approach constants. The normalized solution will converge to $\R^3$.\\

If $T_0=T_2<T_1$, then near $t=T_0$ we have
\begin{equation}
\begin{aligned}
A&\approx k_1\\
B&\approx k_2(T_0-t)^{1/3}\\
C&\approx k_3(T_0-t)^{1/3}\\
D&\approx k_4(T_0-t)^{-1/3}.
\end{aligned}
\end{equation}

Here all sectional curvatures perpendicular to $Y_1$ will have a singularity at $T_0$ of the form $(T_0-t)^{-1}$ while those parallel to $Y_1$ approach constants. The normalized solution will converge to $\R^2$.\\

If $T_0=T_1=T_2$, then near $t=T_0$ we have
\begin{equation}
\begin{aligned}
A&=4(T_0-t)\\
B&\approx k_2(T_0-t)^{1/3}\\
C&\approx k_3(T_0-t)^{1/3}\\
D&\approx k_4(T_0-t)^{-1/3}.
\end{aligned}
\end{equation}

Here all non-zero curvatures will have a singularity of the form $(T_0-t)^{-1}$, and the normalized solution will converge to a manifold $M^2\times\R$ where $M^2$ is the manifold generated by the $2$-forms corresponding to $B$ and $C$.

\subsection{A8. Class $U3I2$.}
Here we may choose a basis for the Lie Albegra $\{X_1,X_2,X_3,X_4\}$ such that the Lie bracket is of the form
\begin{align*}
[X_1,X_2]&=X_3 & [X_1,X_3]&=-X_2 & [X_1,X_4]&=0\\
[X_2,X_3]&=-X_4 & [X_2,X_4]&=0 & [X_3,X_4]&=0.
\end{align*}
This does not correspond to any of the compact homogeneous geometries.\\

Here we diagonalize the metric by letting $Y_i=\Lambda_i^kX_k$ with
\begin{equation*}
\Lambda=\left[\begin{array}{cccc}
1&a_4&a_5&a_6\\0&1&a_2&a_3\\0&0&1&a_1\\0&0&0&1
\end{array}\right].
\end{equation*}
By Proposition 7 in \cite{IJL}, the metric $g(t)$ remains diagonal in the basis $Y_i$ if and only if $a_2=0,~a_1=a_5$ and $a_3=a_4$. We find in this case that $\{Y_i\}$ satisfies the following bracket relations:
\begin{align*}
[Y_1,Y_2]&=Y_3 & [Y_1,Y_3]&=-Y_2 & [Y_1,Y_4]&=0\\
[Y_2,Y_3]&=-Y_4 & [Y_2,Y_4]&=0 & [Y_3,Y_4]&=0.
\end{align*}
Let $g(t)=A(t)\theta_1^2+B(t)\theta_2^2+C(t)\theta_3^2+D(t)\theta_4^2$ where $A(0)=\lambda_1,~B(0)=\lambda_2,~C(0)=\lambda_3$ and $D(0)=\lambda_4$. Then backward Ricci flow reduces to the following system of equations:
\begin{equation}\label{A8eqn}
\begin{aligned}
\frac{dA}{dt}&=-\frac{(B-C)^2}{BC}\\
\frac{dB}{dt}&=\frac{B^2-C^2-AD}{AC}\\
\frac{dC}{dt}&=\frac{C^2-B^2-AD}{AB}\\
\frac{dD}{dt}&=\frac{D^2}{BC}.
\end{aligned}
\end{equation}
By the symmetry of $B$ and $C$ in (\ref{A8eqn}) we may assume that $\lambda_2\geq \lambda_3$. Note also that the equations for $B,C$ and $D$ are identical to those in (\ref{A7ieqn}), so by (\ref{A7iBCD}) and (\ref{A7iD}) we have
\begin{equation}\label{A8BCD}
BCD^2=\lambda_2\lambda_3\lambda_4^2
\end{equation}
and
\begin{equation}\label{A8D}
D=\lambda_4\left(1-\frac{3\lambda_4}{\lambda_2\lambda_3}t\right)^{-1/3}.
\end{equation}

Similar calculations as those used to compute (\ref{A7iADBC}) show that
\begin{equation}\label{A8ADBC}
AD(B+C)=\lambda_1\lambda_4(\lambda_2+\lambda_3).
\end{equation}

Using (\ref{A8eqn}) and (\ref{A8ADBC}), we can solve for $A$:
\begin{equation*}
\frac{dA}{dt}=\frac{-4(k^2+A^2)}{A^2},
\end{equation*}
where
\begin{equation}\label{A8k}
\frac{\lambda_1^2(\lambda_2+\lambda_3)^2}{\lambda_2\lambda_3}=4k^2.
\end{equation}

Now we integrate to find
\begin{equation}\label{A8A}
k\tanh^{-1}\left(\frac{A}{k}\right)-A=-4t+k\tanh^{-1}\left(\frac{\lambda_1}{k}\right)-\lambda_1,
\end{equation}
where $k$ is given in (\ref{A8k}). Near $A=0$ we have
\begin{equation}
A\approx\bigl(12k^2(T_2-t)\bigr)^{1/3}.
\end{equation}

Using (\ref{A8BCD}) and (\ref{A8ADBC}) we find
\begin{align}
B&=\frac{\lambda_4}{2AD}\left(\lambda_1(\lambda_2+\lambda_3)+\sqrt{\lambda_1^2(\lambda_2+\lambda_3)^2-4A^2\lambda_2\lambda_3}\right),\label{A8B}\\
C&=\frac{\lambda_4}{2AD}\left(\lambda_1(\lambda_2+\lambda_3)-\sqrt{\lambda_1^2(\lambda_2+\lambda_3)^2-4A^2\lambda_2\lambda_3}\right),\label{A8C}\\
\intertext{and}
T_0&=\min\left\{\frac{\lambda_2\lambda_3}{3\lambda_4},~\frac{1}{4}\left(k\tanh^{-1}\left(\frac{\lambda_1}{k}\right)-\lambda_1\right)\right\}\notag\\
&=\min{T_1,T_2}.\label{A8T0}
\end{align}

The sectional curvatures are as follows:
\begin{align*}
K(Y_1,Y_2)&=\frac{\frac{B}{C}-3\frac{C}{B}+2}{4A} & K(Y_1,Y_3)&=\frac{\frac{C}{B}-3\frac{B}{C}+2}{4A}\\
K(Y_1,Y_4)&=0 & K(Y_2,Y_3)&=\frac{B^2+C^2-2BC-3AD}{4ABC}\\
K(Y_2,Y_4)&=\frac{D}{4BC} & K(Y_3,Y_4)&=\frac{D}{4BC}.
\end{align*}

If $T_0=T_1<T_2$, then we have the following behavior as $t\rightarrow T_0$:
\begin{equation}
\begin{aligned}
A&\approx k_1\\
B&\approx k_2(T_0-t)^{1/3}\\
C&\approx k_3(T_0-t)^{1/3}\\
D&=k_4(T_0-t)^{-1/3}.
\end{aligned}
\end{equation}
In this case the curvatures parallel to $Y_1$ approach constants while those perpendicular to $Y_1$ will have a singularity of the form $(T_0-t)^{-1}$.\\

If $T_0=T_2<T_1$, then as $t\rightarrow T_0$,
\begin{equation}
\begin{aligned}
A&\approx k_1(T_0-t)^{1/3}\\
B&\approx k_2(T_0-t)^{-1/3}\\
C&\approx k_3(T_0-t)^{1/3}\\
D&\approx k_4.
\end{aligned}
\end{equation}
In this case the curvatures perpendicular to $Y_4$ approach infinity at a rate of $(T_0-t)^{-1}$ while those curvatures parallel to $Y_4$ approach constants.\\

If $T_0=T_1=T_2$, then as $t\rightarrow T_0, AD\rightarrow$ constant, so
\begin{equation}
\begin{aligned}
A&\approx k_1(T_0-t)^{1/3}\\
B&\approx k_2\\
C&\approx k_3(T_0-t)^{2/3}\\
D&= k_4(T_0-t)^{-1/3}.
\end{aligned}
\end{equation}
Here, all non-zero curvatures will have singularities of the form $(T_0-t)^{-1}$.\\

In all of these cases, the volume normalized metric approaches the plane $\R^2$.

\subsection{A9. Class $U3S1$.}

Here we may choose a basis for the Lie Albegra $\{X_1,X_2,X_3,X_4\}$ such that the Lie bracket is of the form
\begin{align*}
[X_1,X_2]&=-X_3 & [X_1,X_3]&=-X_2 & [X_1,X_4]&=0\\
[X_2,X_3]&=X_1 & [X_2,X_4]&=0 & [X_3,X_4]&=0.
\end{align*}

This Lie Algebra structure is a direct sum $\mathfrak{sl}_2\oplus\R$, and the Lie Group structure structure is $(M,G)=(\tilde{SL}(2,\R)\times\R,\tilde{SL}(2,\R)\times\R)$.\\

Here we diagonalize the metric by letting $Y_i=\Lambda_i^kX_k$ with
\begin{equation*}
\Lambda=\left[\begin{array}{cccc}
1&0&0&0\\0&1&0&0\\0&0&1&0\\a_1&a_2&a_3&1
\end{array}\right].
\end{equation*}

\begin{remark}
There are initial metrics in this class which cannot be diagonalized so easy. However, the situations presented in these cases provide us with very complicated situations which we will not address here. See \cite[p.~376-377] {IJL}.
\end{remark}

Here we calculate the curvatures of this diagonalized metric in general form so that we can use it to the specific cases outlined below in Proposition 8. First we calculate the operator $U$ using (\ref{U}):
\begin{align*}
U(Y_i,Y_i)&=0\text{ for all }i\\
U(Y_1,Y_2)&=\frac{B-A}{2C}Y_3+\frac{a_3(B-A)}{2D}Y_4\\
U(Y_1,Y_3)&=\frac{A+C}{2B}Y_2+\frac{A_2(A+C)}{2D}Y_4\\
U(Y_1,Y_4)&=\frac{a_3A}{2B}Y_2-\frac{a_2A}{2C}Y_3\\
U(Y_2,Y_3)&=-\frac{B+C}{2A}Y_1-\frac{a_1(B+C)}{2D}Y_4\\
U(Y_2,Y_4)&=-\frac{a_3B}{2A}Y_1+\frac{a_1B}{2C}Y_3\\
U(Y_3,Y_4)&=-\frac{a_2C}{2A}Y_1+\frac{a_1C}{2B}Y_2.
\end{align*}

The sectional curvatures may then be calculated using (\ref{curvature}):
\begin{equation}\label{A9curvature}
\begin{aligned}
K(Y_i,Y_i)&=0\text{ for all }i\\
K(Y_1,Y_2)&=\frac{1}{4AB}\left[-3C-2B-2A+(A-B)^2\left(\frac{1}{C}+\frac{a_3^2}{D}\right)\right]\\
K(Y_1,Y_3)&=\frac{1}{4AC}\left[-3B-2C+2A+(A+C)^2\left(\frac{1}{B}+\frac{a_2^2}{D}\right)\right]\\
K(Y_1,Y_4)&=\frac{1}{4AD}\left[-3(a_3^2B+a_2^2C)+2A(a_3^2-a_2^2)+A^2\left(\frac{a_3^2}{B}+\frac{a_2^2}{C}\right)\right]\\
K(Y_2,Y_3)&=\frac{1}{4BC}\left[-3A-2C+2B+(B+C)^2\left(\frac{1}{A}+\frac{a_1^2}{D}\right)\right]\\
K(Y_2,Y_4)&=\frac{1}{4BD}\left[-3(a_3^2A+a_1^2C)+2B(a_3^2-a_1^2)+B^2\left(\frac{a_3^2}{A}+\frac{a_1^2}{C}\right)\right]\\
K(Y_3,Y_4)&=\frac{1}{4CD}\left[-3(a_2^2A+a_1^2B)-2C(a_1^2+a_2^2)+C^2\left(\frac{a_2^2}{A}+\frac{a_1^2}{B}\right)\right].
\end{aligned}
\end{equation}

Proposition 8 in \cite{IJL} tells us that when $g_0$ can in fact be diagonalized as above:\\

\begin{tabular}{l}
(i) If $\lambda_1\neq\lambda_2$, the metric remains diagonal if and only if $a_1=a_2=a_3=0$\\
(ii) If $\lambda_1=\lambda_2$, the metric remains diagonal if and only if $a_1=a_2=0$.
\end{tabular}

\subsubsection{A9(i)} Here we consider the case where $\lambda_1\neq\lambda_2$ and $a_1=a_2=a_3=0$. In this case, $\Lambda=I$, and the Lie Bracket relations remain the same. It also happens that this metric is just the product metric $\tilde{SL}(2,\R)\times\R$, so backwards Ricci flow reduces to the case of the three dimensional flow on $\tilde{SL}(2,\R)$. The volume-normalized version of this flow has been discussed in \cite{CSC}.\\

In this case the non-zero curvatures are given by
\begin{equation}\label{A9icurvature}
\begin{aligned}
K(Y_1,Y_2)&=\frac{1}{4ABC}\left[-3C^2-2BC-2AC+A^2-2AB+B^2\right]\\
K(Y_1,Y_3)&=\frac{1}{4ABC}\left[-3B^2-2BC+2AB+A^2+2AC+C^2\right]\\
K(Y_2,Y_3)&=\frac{1}{4ABC}\left[-3A^2-2AC+2AB+B^2+2BC+C^2\right].
\end{aligned}
\end{equation}

Let $Y_i=X_i$, and let $g(t)=A(t)\theta_1^2+B(t)\theta_2^2+C(t)\theta_3^2+D(t)\theta_4^2$ where $A(0)=\lambda_1,~B(0)=\lambda_2,~C(0)=\lambda_3$ and $D(0)=\lambda_4$. Then backward Ricci flow reduces to the following system of equations:
\begin{equation}\label{A9ieqn}
\begin{aligned}
\frac{dA}{dt}&=\frac{A^2-(B+C)^2}{BC}\\
\frac{dB}{dt}&=\frac{B^2-(A+C)^2}{AC}\\
\frac{dC}{dt}&=\frac{C^2-(A-B)^2}{AB}\\
\frac{dD}{dt}&=0.
\end{aligned}
\end{equation}

By the symmetry of the system we may assume that $\lambda_1>\lambda_2$. We also have

\begin{equation*}
\frac{d}{dt}(A-B)=\frac{2}{ABC}(A-B)(A+B+C),
\end{equation*}

so $A>B$ is preserved by (\ref{A9ieqn}).\\

We set
\begin{equation}\label{A9iQ}
Q=\bigl\{(\lambda_1,\lambda_2,\lambda_3)\in\R^3 \bigl| a\geq b>0,~c>0\bigr\}.
\end{equation}
We show that there exists a partition of $Q$ into $Q_1,~Q_2,~S_0$, with $S_0$ a hypersurface in $\R^3$ and $Q_1,~Q_2$ open and connected such that:\\

If $(\lambda_1,\lambda_2,\lambda_3)\in Q_1$, then
\begin{equation}
\begin{aligned}
A(t)&\approx k_1(T_0-t)^{1/3}\\
B(t)&\approx k_2(T_0-t)^{1/3}\\
C(t)&\approx k_3(T_0-t)^{-1/3}\\
D(t)&=\lambda_4.
\end{aligned}
\end{equation}

If $(\lambda_1,\lambda_2,\lambda_3)\in Q_2$, then
\begin{equation}
\begin{aligned}
A(t)&\approx k_1(T_0-t)^{-1/3}\\
B(t)&\approx k_2(T_0-t)^{1/3}\\
C(t)&\approx k_3(T_0-t)^{1/3}\\
D(t)&=\lambda_4.
\end{aligned}
\end{equation}

If $(\lambda_1,\lambda_2,\lambda_3)\in S_0$, then
\begin{equation}
\begin{aligned}
A(t)&\approx k_1\\
B(t)&\approx 4(T_0-t)\\
C(t)&\approx k_1\\
D(t)&=\lambda_4.
\end{aligned}
\end{equation}

In all three cases all the non-zero curvatures have a singularity at $T_0$ of the form $(T_0-t)^{-1}$.\\

We first define
\begin{align}
Q_1&=\bigl\{(\lambda_1,\lambda_2,\lambda_3)\in Q \bigl| C(t_1)\geq A(t_1)\text{ for some time }t_1\geq 0\bigr\}\label{A9iQ1}\\
Q_2&=\bigl\{(\lambda_1,\lambda_2,\lambda_3)\in Q \bigl| C(t_1)\leq A(t_1)-B(t_1)\text{ for some time }t_1\geq 0\bigr\}\label{A9iQ2}\\
S_0&=\bigl\{(\lambda_1,\lambda_2,\lambda_3)\in Q \bigl| A(t)-B(t)\leq C(t)\leq A(t)\text{ for all time }t\geq 0\bigr\}.\label{A9iS0}
\end{align}

Now, if $C(t)\geq A(t)$, then
\begin{equation*}
\frac{d}{dt}\bigl(C-A\bigr)=\frac{1}{ABC}\bigl[(C-A)\bigl((C+A)^2-B^2\bigr)+4ABC\bigr]\geq 4,
\end{equation*}
so $C\geq A$ is preserved.\\

Similarly, if $C(t)\leq A(t)-B(t)$, then
\begin{equation*}
\frac{d}{dt}\bigl(A-B-C\bigr)=\frac{1}{ABC}\bigl[(A-B-C)\bigl(A(A+2B+2C)+(B-C)^2\bigr)+8ABC\bigr]\geq 8,
\end{equation*}
so $C\leq A-B$ is preserved.\\

Thus $Q_1,~Q_2$ and $S_0$ are mutually exclusive sets whose union is all of $Q$. The facts that $Q_1$ and $Q_2$ are open and $S_0$ is a hypersurface in $Q$ are shown in \cite{CGSC}. These results are presented in \cite{CSC} for the normalized Ricci flow, so we give an un-normalized version here to better describe the behavior.\\

We first consider the set $Q_1$, given by (\ref{A9iQ1}). Then $C(t)\geq A(t)$ for all $t\geq t_1$, so we have
\begin{equation*}
\frac{dA}{dt}=\frac{A^2-(B+C)^2}{BC}\leq\frac{-B-2C}{C}<-2,
\end{equation*}
hence $A$ is decreasing, and $T_0<\infty$. Since $B<A<C$, either $\ds\lim_{t\rightarrow T_0}C(t)=\infty$ or $\ds\lim_{t\rightarrow T_0}B(t)=0$. We show that in fact both of these situations happen at the same time. Since $A$ and $B$ are decreasing, then for $C$ large enough we have

\begin{equation}
\frac{d}{dt}(AC)=2(A+B-C)<0.
\end{equation}

Thus $AC$ is bounded above. Thus if $C\rightarrow\infty$ at $T_0$, then it must be the case that $A\rightarrow 0$ at $T_0$, hence also $B\rightarrow 0$ since $B<A$. Now we observe

\begin{equation}
\frac{d}{dt}\bigl(B(C-A)\bigr)=4B>0,
\end{equation}
so $\bigl(B(C-A)\bigr)$ is bounded below. Thus if $B\rightarrow 0$ at $T_0$, then $C-A\rightarrow \infty$, hence also $C\rightarrow\infty$. Thus $A\rightarrow 0,~B\rightarrow 0$ and $C\rightarrow\infty$ at time $t=T_0$. Now,

\begin{equation}
\frac{d}{dt}(AC-AB)=4A,
\end{equation}
which is positive and approaches $0$ at $t=T_0<\infty$. Therefore we know that $AC-AB\rightarrow k_A>0$ as $t\rightarrow T_0$. Since $AB\rightarrow 0$, we see that
\begin{equation}\label{A9ikA}
AC\rightarrow k_A.
\end{equation}

Similarly, as $t\rightarrow T_0$,
\begin{equation}\label{A9ikB}
BC\rightarrow k_B>0.
\end{equation}

Now by (\ref{A9ieqn}) we see that
\begin{equation}
\frac{dC}{dt}\sim\frac{C^2}{AB}\sim kC^4,
\end{equation}
hence by the same argument used to solve (\ref{A3A1}), we see
\begin{equation}\label{A9iC}
C\approx k_3(T_0-t)^{-1/3}.
\end{equation}

Now combining (\ref{A9ikA}), (\ref{A9ikB}) and (\ref{A9iC}), we have the solution to (\ref{A9ieqn}):
\begin{equation}\label{A9iSol1}
\begin{aligned}
A(t)&\approx k_1(T_0-t)^{1/3}\\
B(t)&\approx k_2(T_0-t)^{1/3}\\
C(t)&\approx k_3(T_0-t)^{-1/3}\\
D(t)&=\lambda_4.
\end{aligned}
\end{equation}

The normalized solution approaches the plane $\R^2$.\\

Now we consider the set $Q_2$ given by (\ref{A9iQ2}). Here there is some time $t_1$ such that for $t\geq t_1\geq 0,~ (A-B)<C$. We denote $A_1=A(t_1),~B_1=B(t_1)$ and $C_1=C(t_1)$. We see that for $t\geq t_1$
\begin{align*}
\frac{dC}{dt}&=\frac{C^2-(A+B)^2}{AB}<0,
\intertext{and}
\frac{dA}{dt}&=\frac{A^2-(B-C)^2}{BC}>0,
\end{align*}
so $C$ is decreasing and $A$ is increasing. Thus, as $t\rightarrow T_0$, either $A\rightarrow\infty,~B\rightarrow 0$ or $C\rightarrow 0$. We calculate
\begin{equation}
\frac{d}{dt}\bigl(A(B+C)\bigr)=-4B-4C<0,
\end{equation}
so $\bigl(A(B+C)\bigr)$ is bounded above. Thus if $A\rightarrow\infty$ as $t\rightarrow T_0$, then we know $(B+C)\rightarrow 0$ at $t=T_0$, hence $B\rightarrow 0$ and $C\rightarrow 0$. Now for $t\geq t_1$,
\begin{equation}
\frac{d}{dt}\left(\frac{B}{C}\right)=\frac{2}{A}(B+C)(B-C-A)<0,
\end{equation}
so $\ds\frac{B}{C}$ is bounded above and if $C\rightarrow 0$ then $B\rightarrow 0$. Lastly,
\begin{equation}\label{A9iBABC}
\frac{d}{dt}\bigl(B(A-B-C)\bigr)=\frac{2B}{AC}(A^2+C^2-B^2)>0,
\end{equation}
so $\bigl(B(A-B-C\bigr)$ is bounded below, and if $B\rightarrow 0$ then $(A-B-C)\rightarrow \infty$, hence also $A\rightarrow\infty$. Thus at $t=T_0$, $A\rightarrow\infty,~B\rightarrow 0$ and $C\rightarrow 0$.\\

Now by (\ref{A9iBABC}) we see that $(AB-BC-B^2)$ is increasing for $t\geq t_1$.Thus we see that $AB>B_1(A_1-B_1-C_1)$. However,
\begin{equation}
\frac{d}{dt}(AB)=-2(A+B+C)<0,
\end{equation}
so $AB$ is decreasing and bounded below by a positive number. Thus at $t\rightarrow T_0$, we have
\begin{equation}\label{A9ikAB}
AB\rightarrow k_{AB}
\end{equation}
for some positive constant $k_{AB}$. Now
\begin{equation}\label{A9iAC}
\frac{d}{dt}(AC)=2(A-B-C)>0,
\end{equation}
so $(AC)$ is increasing and bounded below. Thus, since $B$ and $C$ are approaching $0$ and $A$ is approaching infinity at $T_0$, then by (\ref{A9ieqn}) and (\ref{A9ikAB}) there is some positive number $k$ such that
\begin{equation}
\frac{dA}{dt}=\frac{A^2-(B+C)^2}{BC}\leq kA^4,
\end{equation}
hence we know
\begin{equation}
A< \left(A_1^{-3}-3kt\right)^{-1/3}.
\end{equation}
By (\ref{A9iAC}) we see that
\begin{equation}
\frac{d}{dt}(AC)<2A< 2\bigl(A_1^{-3}-3kt\bigr)^{-1/3}.
\end{equation}

Upon integrating we discover $AC<\infty$, and since $AC$ is increasing, then we see that at $t\rightarrow T_0$, we have
\begin{equation}\label{A9ikAC}
AC\rightarrow k_{AC}
\end{equation}
for some positive constant $k_{AC}$. Now by (\ref{A9ieqn}), (\ref{A9ikAB}) and (\ref{A9ikAC}),
\begin{equation}\label{A9iA}
A\sim kA^4.
\end{equation}

Solving (\ref{A9iA}) the same way we solved (\ref{A3A1}), and then using (\ref{A9ikAB}) and (\ref{A9ikAC}), we have the following end behavior of the solutions to (\ref{A9ieqn}):
\begin{equation}
\begin{aligned}
A(t)&\approx k_1(T_0-t)^{-1/3}\\
B(t)&\approx k_2(T_0-t)^{1/3}\\
C(t)&\approx k_3(T_0-t)^{1/3}\\
D(t)&=\lambda_4.
\end{aligned}
\end{equation}

The normalized solution approaches the plane $\R^2$.\\

Now we consider the set $S_0$ given by (\ref{A9iS0}). Here we know 
\begin{equation*}
A(t)-B(t)<C(t)<A(t)
\end{equation*}
for all time $0\leq t< T_0$. Thus we have

\begin{align*}
\frac{dA}{dt}&=\frac{A^2-(B+C)^2}{BC}<0,
\intertext{and}
\frac{dC}{dt}&=\frac{C^2-(A-B)^2}{AB}>0,
\end{align*}
so $A$ is decreasing while $C$ is increasing.

However, $A(t)>C(t)$, so both $A$ and $C$ are approaching constants. Now
\begin{equation*}
\frac{dB}{dt}=\frac{B^2-(A+C)^2}{AC}<-2,
\end{equation*}
so $B$ approaches $0$ in finite time $T_0<\infty$. Since $A-B<C<A$, then in fact $A$ and $C$ approach the same constant, $k_1$, at time $t=T_0$. Thus near $t=T_0$ we have the approximation:
\begin{equation}
\frac{dB}{dt}\approx -\frac{(A+C)^2}{AC}\approx -4.
\end{equation}

Thus near $t=T_0$ we have the following behavior:
\begin{equation}
\begin{aligned}
A&\approx k_1\\
B&\approx -4(t-T_0)\\
C&\approx k_1\\
D&=\lambda_4.
\end{aligned}
\end{equation}

The normalized solution approaches $\R^3$.

\subsubsection{A9(ii)}

We find in this case that $\{Y_i\}$ satisfies the following bracket relations:
\begin{center}
\begin{tabular}{lll}
$[Y_1,Y_2]=-Y_3$\hspace{.2in} & $[Y_1,Y_3]=-Y_2$\hspace{.2in} & $[Y_1,Y_4]=-a_3Y_2$\\
$[Y_2,Y_3]=Y_1$\hspace{.2in} & $[Y_2,Y_4]=a_3Y_1$\hspace{.2in} & $[Y_3,Y_4]=0$
\end{tabular}
\end{center}
Let $g(t)=A(t)\theta_1^2+B(t)\theta_2^2+C(t)\theta_3^2+D(t)\theta_4^2$ where $A(0)=\lambda_1,~B(0)=\lambda_2,~C(0)=\lambda_3$ and $D(0)=\lambda_4$. Then backward Ricci flow reduces to the following system of equations:
\begin{equation}\label{A9iieqn1}
\begin{aligned}
\frac{dA}{dt}&=\frac{A^2-(B+C)^2}{BC}+\frac{A^2-B^2}{BD}a_3^2\\
\frac{dB}{dt}&=\frac{B^2-(A+C)^2}{AC}+\frac{B^2-A^2}{AD}a_3^2\\
\frac{dC}{dt}&=\frac{C^2-(A-B)^2}{AB}\\
\frac{dD}{dt}&=-\frac{(A+B)^2}{AB}a_3^2\\
\lambda_1&=\lambda_2.
\end{aligned}
\end{equation}

We calculate
\begin{align*}
\frac{d}{dt}\bigl(A-B\bigr)&=\biggl[\frac{1}{ABC}\bigl(-C^2+A^2+2AB+B^2\bigr)\\
&\qquad+\frac{a_3^2}{ABD}\bigl(A^2+2AB+B^2\bigr)\biggr]\bigl(A-B\bigr),
\end{align*}
so $A=B$ for all $t$.\\

Now (\ref{A9iieqn1}) reduces to
\begin{equation}\label{A9iieqn}
\begin{aligned}
\frac{dA}{dt}&=-\frac{C}{A}-2\\
B&=A\\
\frac{dC}{dt}&=\frac{C^2}{A^2}\\
\frac{dD}{dt}&=-4a_3^2.
\end{aligned}
\end{equation}

It is clear that $D=\lambda_4-4a_3^2t$, so $\ds T_0\leq\frac{\lambda_4}{4a_3^2}$. Also, letting $f'$ denote $\ds\frac{df}{dt}$ we may compute
\begin{equation*}
\frac{A'}{(A'+1)(A'+2)}A''=-\frac{2}{A}A'.
\end{equation*}
Integrating gives
\begin{equation}\label{A9iiAprime}
A'=\frac{-\Lambda-4A^2-\sqrt{\Lambda^2+4\Lambda A^2}}{2A^2}.
\end{equation}
where
\begin{equation}\label{A9iiLambda}
\Lambda=\frac{\lambda_3^2\lambda_1}{\lambda_3+\lambda_1}.
\end{equation}

Now we may calculate
\begin{equation}\label{A9iiA}
\frac{A}{2}-\frac{\sqrt{\Lambda}}{4}\sinh^{-1}\left(\frac{2A}{\sqrt{\Lambda}}\right)=T_2-t,
\end{equation}
so near $t=T_2$,
\begin{equation}
A\approx\left(3e^{\Lambda}(T_2-t)\right)^{1/3}.
\end{equation}

By (\ref{A9iieqn}) we know that near $A=0$,
\begin{equation}
C=-A\left(\frac{dA}{dt}+2\right)\approx \frac{\Lambda}{2}\left(3e^{\Lambda}(T_2-t)\right)^{-1/3}.
\end{equation}

We observe that
\begin{align}
T_0&=\min\left\{\frac{\lambda_4}{4a_3^2},~\frac{\lambda_1}{2}-\frac{\lambda_3}{4}\sqrt{\frac{\lambda_1}{\lambda_1+\lambda_3}\sinh^{-1}\biggl(\frac{2}{\lambda_3}\sqrt{\lambda_1(\lambda_1+\lambda_3)}\biggr)}\right\}\\
&=\min\bigl\{T_1,T_2\bigr\}.
\end{align}

We calculate the sectional curvatures using (\ref{A9curvature}) with $a_1=a_2=0$ and $A=B$:
\begin{align*}
K(Y_1,Y_2)&=-\frac{4A+3C}{4A^2} & K(Y_1,Y_3)&=\frac{C}{4A^2} & K(Y_1,Y_4)&=0\\
K(Y_2,Y_3)&=\frac{C}{4A^2} & K(Y_2,Y_4)&=0 & K(Y_3,Y_4)&=0.
\end{align*}

If $T_0=T_1<T_2$, then near $t=T_0$ we have
\begin{equation}
\begin{aligned}
A=B&\approx k_1\\
C&\approx k_3\\
D&=k_4(T_0-t).
\end{aligned}
\end{equation}

Here all sectional curvatures will also approach constants as $t\rightarrow T_0$, and the volume normalized solution approaches the hyperplane $\R^3$.\\

If $T_0=T_2< T_1$, then near $t=T_0$ we have
\begin{equation}
\begin{aligned}
A=B&\approx k_1(T_0-t)^{1/3}\\
C&\approx k_3(T_0-t)^{-1/3}\\
D&\approx k_4.
\end{aligned}
\end{equation}

Here all non-zero sectional curvatures will have a singularity of the form $(T_0-t)^{-1}$. The volume normalized solution approaches $\R^2$.\\

If $T_0=T_1= T_2$, then near $t=T_0$ we have
\begin{equation}
\begin{aligned}
A=B&\approx k_1(T_0-t)^{1/3}\\
C&\approx k_3(T_0-t)^{-1/3}\\
D&=k_4(T_0-t).
\end{aligned}
\end{equation}

Here all non-zero sectional curvatures will have a singularity of the form $(T_0-t)^{-1}$, and the volume normalized solution approaches the plane $\R\times M$, where $M$ is a $2$-manifold generated by the $2$-forms corresponding to $A$ and $B$.

\subsection{A10. Class $U3S3$.}
Here we may choose a basis for the Lie Albegra $\{X_1,X_2,X_3,X_4\}$ such that the Lie bracket is of the form
\begin{align*}
[X_1,X_2]&=X_3 & [X_1,X_3]&=-X_2 & [X_1,X_4]&=0\\
[X_2,X_3]&=X_1 & [X_2,X_4]&=0 & [X_3,X_4]&=0.
\end{align*}

This Lie Algebra structure is a direct sum $\mathfrak{su}(2)\oplus\R$, and the Lie Group structure structure is $(M,G)=(S^3\times\R,SU(2)\times\R)$.\\

Here we diagonalize the metric by letting $Y_i=\Lambda_i^kX_k$ with
\begin{equation*}
\Lambda=\left[\begin{array}{cccc}
1&0&0&0\\0&1&0&0\\0&0&1&0\\a_1&a_2&a_3&1
\end{array}\right].
\end{equation*}

Here I calculate the curvatures of this diagonalized metric in general form so that I can use it to the specific cases outlined below in Proposition 9. First I calculate the operator $U$ using (\ref{U}):
\begin{align*}
U(Y_i,Y_i)&=0\text{ for all }i\\
U(Y_1,Y_2)&=\frac{B-A}{2C}Y_3+\frac{a_3(B-A)}{2D}Y_4\\
U(Y_1,Y_3)&=\frac{A-C}{2B}Y_2+\frac{A_2(A-C)}{2D}Y_4\\
U(Y_1,Y_4)&=\frac{a_3A}{2B}Y_2-\frac{a_2A}{2C}Y_3\\
U(Y_2,Y_3)&=\frac{C-B}{2A}Y_1+\frac{a_1(C-B)}{2D}Y_4\\
U(Y_2,Y_4)&=-\frac{a_3B}{2A}Y_1+\frac{a_1B}{2C}Y_3\\
U(Y_3,Y_4)&=\frac{a_2C}{2A}Y_1-\frac{a_1C}{2B}Y_2.
\end{align*}

The sectional curvatures may then be calculated using (\ref{curvature}):
\begin{equation}\label{A10curvature}
\begin{aligned}
K(Y_i,Y_i)&=0\text{ for all }i\\
K(Y_1,Y_2)&=\frac{1}{4AB}\left[-3C+2B+2A+(A-B)^2\left(\frac{1}{C}+\frac{a_3^2}{D}\right)\right]\\
K(Y_1,Y_3)&=\frac{1}{4AC}\left[-3B+2C+2A+(A-C)^2\left(\frac{1}{B}+\frac{a_2^2}{D}\right)\right]\\
K(Y_1,Y_4)&=\frac{1}{4AD}\left[-3(a_3^2B+a_2^2C)+2A(a_2^2+a_3^2)+A^2\left(\frac{a_3^2}{B}+\frac{a_2^2}{C}\right)\right]\\
K(Y_2,Y_3)&=\frac{1}{4BC}\left[-3A+2C+2B+(B-C)^2\left(\frac{1}{A}+\frac{a_1^2}{D}\right)\right]\\
K(Y_2,Y_4)&=\frac{1}{4BD}\left[-3(a_3^2A+a_1^2C)+2B(a_3^2+a_1^2)+B^2\left(\frac{a_3^2}{A}+\frac{a_1^2}{C}\right)\right]\\
K(Y_3,Y_4)&=\frac{1}{4CD}\left[-3(a_2^2A+a_1^2B)+2C(a_1^2+a_2^2)+C^2\left(\frac{a_2^2}{A}+\frac{a_1^2}{B}\right)\right].
\end{aligned}
\end{equation}
Proposition 9 in \cite{IJL} breaks this down into three situations:\\

(i) If $\lambda_1=\lambda_2=\lambda_3$, them the metric remains diagonal for all choices of $a_1,~a_2$ and $a_3$.\\

(ii) If $\lambda_i\neq\lambda_j=\lambda_k$ for some permutation $\{i,j,k\}$ of $\{1,2,3\}$, then the metric remains diagonal if and only if $a_j=a_k=0$.\\

(iii) If $\lambda_1,~\lambda_2$ and $\lambda_3$ are all different, then the metric will remain diagonal if and only if $a_1=a_2=a_3=0$.\\

In all three of these cases,backward Ricci flow reduces to the system:
\begin{equation}\label{A10eqn}
\begin{aligned}
\frac{dA}{dt}&=\frac{A^2-(B-C)^2}{BC}\\
\frac{dB}{dt}&=\frac{B^2-(A-C)^2}{AC}\\
\frac{dC}{dt}&=\frac{C^2-(A-B)^2}{AB}\\
\frac{dD}{dt}&=0\\
\end{aligned}
\end{equation}

Thus this reduces to $3$ dimensions, and the normalized flow is analyzed in \cite{CSC}.\\

By symmetry, we may in fact assume $\lambda_1\geq\lambda_2\geq\lambda_3$. Now
\begin{equation}
\frac{d}{dt}(A-B)=\frac{1}{ABC}(A-B)\bigl((A+B)^2-C^2\bigr),
\end{equation}
so the conditions $A(t)\geq B(t)$ and $A(t)=B(t)$ are preserved under backwards Ricci flow. Similarly, the conditions $B(t)\geq C(t)$ and $B(t)=C(t)$ are preserved.
\subsubsection{A10(i)}
If $\lambda_1=\lambda_2=\lambda_3$, then $A\equiv B\equiv C$. Then (\ref{A10eqn}) reduces to
\begin{equation}
\begin{aligned}
\frac{dA}{dt}=\frac{dB}{dt}=\frac{dC}{dt}&=1\\
\frac{dD}{dt}&=0,
\end{aligned}
\end{equation}

and the solution, which exists for all time, is
\begin{equation}\label{A10Sol1}
\begin{aligned}
A(t)=B(t)=C(t)&=\lambda_1+t,\\
D(t)&=\lambda_4.
\end{aligned}
\end{equation}

Using (\ref{A10curvature}) with $A=B=C$ we find the non-zero curvatures are
\begin{equation}
K(Y_1,Y_2)=K(Y_1,Y_3)=K(Y_2,Y_3)=\frac{1}{4A}.
\end{equation}

Thus all non-zero curvatures approach $0$ at a rate of $4t^{-1}$.\\

The solution here is actually just a product of an expanding $3$-sphere with a line. The volume normalized solution converges to the hyperplane $\R^3$.

\subsubsection{A10(ii)}
If $\lambda_1=\lambda_2>\lambda_3$, then $A(t)=B(t)> C(t)$ for all $0<t<T_0$. Then (\ref{A10eqn}) reduces to 
\begin{equation}\label{A10eqnii}
\begin{aligned}
\frac{dA}{dt}&=2-\frac{C}{A}\\
B&=A\\
\frac{dC}{dt}&=\frac{C^2}{A^2}\\
\frac{dD}{dt}&=0.
\end{aligned}
\end{equation}

Denoting $\frac{df}{dt}$ by $f'$ we calculate

\begin{equation*}
\frac{-A'}{(A'-1)(A'-2)}A''=\frac{2}{A}A'.
\end{equation*}

Using the fact that $A'>1$, (\ref{A10a}) we calculate
\begin{equation}\label{A10Aprime}
A'=\frac{1+4\Lambda A^2-\sqrt{1+4\Lambda A^2}}{2\Lambda A^2}.
\end{equation}
where
\begin{equation}\label{A10Lambda}
\Lambda=\frac{\lambda_1-\lambda_3}{\lambda_1\lambda_3^2}.
\end{equation}
Upon integrating we have
\begin{equation}\label{A10A}
\frac{1}{2}A+\frac{1}{4\sqrt{\Lambda}}\sinh^{-1}\bigl(2\sqrt{\Lambda}A\bigr)=t+\frac{1}{2}+\frac{1}{4\sqrt{\Lambda}}\sinh^{-1}\bigl(2\sqrt{\Lambda}\lambda_1\bigr),
\end{equation}
so as $A$ approaches infinity, $t\approx \frac{1}{2}A$.\\

By (\ref{A10eqnii}), $\ds A'=2-\frac{C}{A}$  so
\begin{equation}\label{A10AtoC}
C=\frac{-1+\sqrt{1+4\Lambda A^2}}{2\Lambda A}.
\end{equation}

Now by equations (\ref{A10eqnii}), (\ref{A10A}) and (\ref{A10AtoC}), we have the following behavior as $t\rightarrow\infty$:
\begin{equation}\label{A10Sol2}
\begin{aligned}
A(t)&=B(t)\approx 2t\\
C(t)&\rightarrow\frac{1}{\sqrt{\Lambda}}\\
D(t)&=\lambda_4.
\end{aligned}
\end{equation}
where $\Lambda$ is given in (\ref{A10Lambda}).\\

Using (\ref{A10curvature}) with $a_1=a_2=0$ and $A=B$ we find the non-zero curvatures are given by
\begin{equation}
\begin{aligned}
K(Y_1,Y_2)&=\frac{4A-3C}{4A^2}\\
K(Y_1,Y_3)&=\frac{C}{4A^2}\\
K(Y_2,Y_3)&=\frac{C}{4A^2}.
\end{aligned}
\end{equation}

Thus all non-zero curvatures parallel to $Y_3$ approach $0$ at a rate of $t^{-2}$ while the remaining nonzero curvatures approach $0$ at a rate of $t^{-1}$. The normalized solution will converge to the plane $\R^2$.\\

If $\lambda_1>\lambda_2=\lambda_3$, then $A(t)>B(t)=C(t)$ for all $0<t<T_0$, and (\ref{A10eqn}) reduces to the following:
\begin{equation}\label{A10eqn3}
\begin{aligned}
\frac{dA}{dt}&=\frac{A^2}{B^2}\\
\frac{dB}{dt}&=2-\frac{A}{B}\\
C&=B\\
\frac{dD}{dt}&=0.
\end{aligned}
\end{equation}

Similarly to the case $A=B>C$ we calculate

\begin{align}
B''&=-\frac{2}{B}(B'-1)(B'-2),\notag\\
B'&=\frac{4\Lambda B^2-1-\sqrt{1-4\Lambda B^2}}{2\Lambda B^2}.\label{A10Bprime}
\end{align}
where
\begin{equation}\label{A10Lambdab}
\Lambda=\frac{\lambda_1-\lambda_2}{\lambda_1^2\lambda_2}.
\end{equation}
Upon integrating, and using the fact that $B'<1$ we have
\begin{equation}\label{A10B}
\frac{1}{2}B-\frac{1}{4\sqrt{\Lambda}}\sin^{-1}\bigl(2\sqrt{\Lambda}B\bigr)=t-T_0,
\end{equation}
where
\begin{equation}\label{A10T0}
T_0=\frac{\lambda_2}{2}-\frac{\lambda_1\sqrt{\lambda_2}}{4\sqrt{\lambda_1-\lambda_2}}\sin^{-1}\left(\frac{2\sqrt{\lambda_2(\lambda_1-\lambda_2)}}{\lambda_1}\right).
\end{equation}

As $B$ approaches $0$, we observe that $(t-T_0)\approx-\frac{\Lambda}{3}B^3$.

Note here that since $\ds B'=2-\frac{A}{B}$, then
\begin{equation}\label{A10BtoA}
A=\frac{1+\sqrt{1-4\Lambda B^2}}{2\Lambda B}.
\end{equation}

Thus by equations (\ref{A10eqn3}) and (\ref{A10BtoA}), we have the following behavior as $t\rightarrow T_0$:
\begin{equation}
\begin{aligned}
A(t)&\approx\frac{1}{\Lambda}\left(\frac{3}{\Lambda}(T_0-t)\right)^{-1/3}\\
B(t)&=C(t)\approx \left(\frac{3}{\Lambda}(T_0-t)\right)^{1/3}\\
D(t)&=\lambda_4.
\end{aligned}
\end{equation}
where $\Lambda$ is given in (\ref{A10Lambdab}), and $T_0$ is given in (\ref{A10T0}).\\

With $a_2=a_3=0$ and $B=C$, (\ref{A10curvature}) tells us the non-zero curvatures are given by
\begin{equation}
\begin{aligned}
K(Y_1,Y_2)&=\frac{A}{4B^2}\\
K(Y_1,Y_3)&=\frac{A}{4B^2}\\
K(Y_2,Y_3)&=\frac{4B-3A}{4B^2}.
\end{aligned}
\end{equation}

Thus all non-zero curvatures are of the form $(T_0-t)^{-1}$. The volume normalized solution approaches the plane $\R^2$.
\subsubsection{A10(iii)}

In the case $\lambda_1>\lambda_2>\lambda_3$, we have by \cite{CSC} that the end behavior near $t=T_0$ is the same as when $\lambda_1>\lambda_2=\lambda_3$. However, we do not have more explicit solutions like we do in the special case. Since the solutions in \cite{CSC} are using normalized backward Ricci flow, I shall present a slightly different argument here.\\

We know that for all $t$ we have $A(t)>B(t)>C(t)$. From (\ref{A10eqn}) we calculate

\begin{equation}\label{A10Clinear}
\frac{dC}{dt}=\frac{C^2-A^2+2AB-B^2}{AB}<2-\frac{A}{B}<1,
\end{equation}

so $C(t)<\lambda_3+t$. Now we may calculate

\begin{align*}
\frac{d}{dt}\left(\frac{A-B}{C}\right)&=\frac{2}{ABC}(A-B)(A^2+B^2-C^2)\\
&>\frac{2A^2(A-B)}{ABC^2}\\
&>2\left(\frac{A-B}{C}\right)\cdot\frac{1}{C}\\
&>2\left(\frac{A-B}{C}\right)\cdot\frac{1}{t+\lambda_3},
\end{align*}
so
\begin{equation*}
\ln\left(\frac{A-B}{C}\right)>2\ln(t+\lambda_3)+\ln\left(\frac{\lambda_1-\lambda_2}{\lambda_3}\right),
\end{equation*}
and we have that
\begin{equation*}
\left(\frac{A-B}{C}\right)>\left(\frac{\lambda_1-\lambda_2}{\lambda_3}\right)(t+\lambda_3)^2.
\end{equation*}

Therefore $\ds\frac{A-B}{C}$ increases at least quadratically. Thus either $T_0<\infty$ or $\ds\frac{dC}{dt}=\frac{C^2-(A-B)^2}{AB}<0$ for all $t$ large enough.\\

Either way, we know by (\ref{A10Clinear}) that $C$ is bounded, hence $C\leq k_C$ for some $0<k_C<\infty$. Now

\begin{align}
\frac{d}{dt}\left(\frac{A}{B}\right)&=\frac{2(A-B)}{B^2C}(A+B-C)\notag\\
&>\frac{2A(A-B)}{k_CB^2}\notag\\
&=\frac{2}{k_C}\left(\frac{A}{B}\right)\left(\frac{A}{B}-1\right).\label{A10AoverB}
\end{align}

Thus we see that $\ds\frac{A}{B}\rightarrow\infty$ in finite time, so we know $T_0<\infty$. Since for all $t$ we know $A>B>C$, then near $t=T_0$ either $A\rightarrow\infty$ or $C\rightarrow 0$. Now we may calculate

\begin{equation*}
\frac{d}{dt}\left(\frac{B}{C}\right)=\frac{2}{AC^2}(B-C)(B+C-A),
\end{equation*}
which is negative for $t$ close enough to $T_0$. Thus we know that
\begin{equation}\label{A10BClimit}
\lim_{t\rightarrow T_0}\frac{B}{C}=k_{BC}
\end{equation}
for some $k_{BC}\geq 1$. Thus if $C\rightarrow 0$ then $B\rightarrow 0$ as well. Now we observe
\begin{equation}
\frac{d}{dt}(ABC)=-(A^2+B^2+C^2)+2(AB+AC+BC),
\end{equation}
which is negative for $A$ large enough and $B$ and $C$ bounded. Thus if $A\rightarrow\infty$ as $t\rightarrow T_0$, then $ABC$ is bounded, hence $C\rightarrow 0$. Similarly,
\begin{equation}\label{A10AC}
\frac{d}{dt}(AC)=2(A+C-B),
\end{equation}
which is positive for $B$ close enough to $0$ and $A$ bounded below. Hence if $B\rightarrow 0$ then also $A\rightarrow\infty$. Thus we know that as $t\rightarrow T_0$ we have $A\rightarrow\infty,~B\rightarrow 0$ and $C\rightarrow 0$.\\

Now by (\ref{A10eqn}) we observe
\begin{equation}\label{A10A1}
\frac{dA}{dt}=\frac{A^2-(B-C)^2}{BC}<\frac{A^2}{BC}.
\end{equation}

Also, by (\ref{A10AC}) we know that $AC$ is bounded from below. Similarly, $AB$ is bounded from below, so (\ref{A10A1}) becomes
\begin{equation*}
\frac{dA}{dt}<kA^4,
\end{equation*}
hence
\begin{equation}\label{A10A2}
A(t)<\bigl(\lambda_1^{-3}-3kt\bigr)^{-1/3}.
\end{equation}

Now we have
\begin{equation}\label{A10AB}
\frac{d}{dt}(AB)=2(A+B-C)<\tilde{k}A<\tilde{k}\bigl(\lambda_1^{-3}-3kt\bigr)^{-1/3}.
\end{equation}
Integrating (\ref{A10AB}) gives us
\begin{equation*}
AB<\lambda_1\lambda_2+\frac{\tilde{k}}{2\lambda_1^2}-\frac{\tilde{k}}{2}\bigl(\lambda_1^{-3}-3kt\bigr)^{2/3}.
\end{equation*}
hence
\begin{equation}\label{A10ABlimit}
\lim_{t\rightarrow T_0}(AB)=k_{AB}.
\end{equation}
for some $0<k_{AB}<\infty$. Combining (\ref{A10BClimit}) and (\ref{A10ABlimit}) we conclude
\begin{equation}\label{A10AClimit}
\lim_{t\rightarrow T_0}(AC)=\frac{k_{AB}}{k_{BC}}.
\end{equation}

Thus we see that near $t=T_0$, 
\begin{equation*}
BC\sim\frac{k_{AB}^2}{k_{BC}A^2},
\end{equation*}
hence (\ref{A10eqn}) tells us
\begin{equation}\label{A10A2}
\frac{dA}{dt}\sim\frac{k_{BC}}{k_{AB}^2}A^4.
\end{equation}
Solving as we did for equation $(\ref{A3A1})$ gives
\begin{equation}\label{A10A3}
A\approx k_1(T_0-t)^{-1/3}.
\end{equation}

Now using (\ref{A10ABlimit}) and (\ref{A10AClimit}) we have the end behavior of (\ref{A10eqn}) near $t=T_0$:
\begin{equation}\label{A10Sol}
\begin{aligned}
A(t)&\approx k_1(T_0-t)^{-1/3}\\
B(t)&\approx k_2(T_0-t)^{1/3}\\
C(t)&\approx k_3(T_0-t)^{1/3}\\
D(t)&=\lambda_4.
\end{aligned}
\end{equation}

From (\ref{A10curvature}) with $a_1=a_2=a_3=0$ we have that the non-zero curvatures are
\begin{equation}
\begin{aligned}
K(Y_1,Y_2)&=\frac{1}{4ABC}(-3C^2+2AC+2BC+A^2-2AB+B^2)\\
K(Y_1,Y_3)&=\frac{1}{4ABC}(-3B^2+2AB+2BC+A^2-2AC+C^2)\\
K(Y_2,Y_3)&=\frac{1}{4ABC}(-3A^2+2AB+2BC+B^2-2BC+C^2).
\end{aligned}
\end{equation}

Thus all the non-zero curvatures have a singularity of the form $(T_0-t)^{-1}$, and the volume-normalized solution approaches $\R^2$.

\section{The non-Bianchi cases}

In this section we examine the compact locally homogeneous geometries whose isotropy group is not trivial, so the dimension of the Lie Group is higher than the dimension of the manifold. Again, these cases are well-understood, but I include them here for completion. The following can be found in \cite{IJL}.

\subsection{B1. $H^3\times\R$}
Any initial metric can be written as
\begin{equation}
g_0=R^2g_{H^3}+du^2
\end{equation}
for some $R>0$. The solution to backward Ricci flow is given by 
\begin{equation}
g(t)=(R^2-4t)g_{H^3}+du^2,\hspace{.2in}-\infty<t<\frac{R^2}{4}.
\end{equation}
\subsection{B2. $S^2\times\R^2$}
Any initial metric can be written as
\begin{equation}
g_0=R^2g_{S^3}+du_1^2+du_2^2
\end{equation}
for some $R>0$. The solution to backward Ricci flow is given by 
\begin{equation}
g(t)=(R^2+2t)g_{S^2}+du_1^2+du_2^2,\hspace{.2in}-\frac{R^2}{2}<t<\infty.
\end{equation}
\subsection{B3. $H^2\times\R^2$}
Any initial metric can be written as
\begin{equation}
g_0=R^2g_{H^2}+du_1^2+du_2^2
\end{equation}
for some $R>0$. The solution to backward Ricci flow is given by 
\begin{equation}
g(t)=(R^2-2t)g_{H^2}+du_1^2+du_2^2,\hspace{.2in}-\infty<t<\frac{R^2}{2}.
\end{equation}
\subsection{B4. $S^2\times S^2$}
Any initial metric can be written as
\begin{equation}
g_0=R_1^2g_{S^2}+R_2^2g_{S^2}
\end{equation}
for some $R_1>0,~R_2>0$. The solution to backward Ricci flow is given by 
\begin{equation}
g(t)=(R_1^2+2t)g_{S^2}+(R_2^2+2t)g_{S^2},\hspace{.2in}-\min\left\{\frac{R_1^2}{2},\frac{R_2^2}{2}\right\}<t<\infty.
\end{equation}
\subsection{B5. $S^2\times H^2$}
Any initial metric can be written as
\begin{equation}
g_0=R_1^2g_{S^2}+R_2^2g_{H^2}
\end{equation}
for some $R_1>0,~R_2>0$. The solution to backward Ricci flow is given by 
\begin{equation}
g(t)=(R_1^2-2t)g_{S^2}+(R_2^2+2t)g_{H^2},\hspace{.2in}-\frac{R_1^2}{2}<t<\frac{R_2^2}{2}.
\end{equation}
\subsection{B6. $H^2\times H^2$}
Any initial metric can be written as
\begin{equation}
g_0=R_1^2g_{H^2}+R_2^2g_{H^2}
\end{equation}
for some $R_1,R_2>0,~R_1\neq R_2$. The solution to backward Ricci flow is given by 
\begin{equation}
g(t)=(R_1^2+2t)g_{H^2}+(R_2^2+2t)g_{H^2},\hspace{.2in}-\infty<t<\min\left\{\frac{R_1^2}{2},\frac{R_2^2}{2}\right\}.
\end{equation}
\subsection{B7. $\C P^2$}
Any initial metric can be written as
\begin{equation}
g_0=R^2g_{FS}
\end{equation}
for some $R>0$, where $g_{FS}$ is the Fubini-Study metric on complex projective space, $\C P^2$, with constant holomorphic bisectional curvature $1$. The solution to backward Ricci flow is given by 
\begin{equation}
g(t)=(R^2+6t)g_{FS},\hspace{.2in}-\frac{R^2}{6}<t<\infty.
\end{equation}
\subsection{B8. $\C H^2$}
Any initial metric can be written as
\begin{equation}
g_0=R^2g_{\C H^2}
\end{equation}
for some $R>0$, where $g_{\C H^2}$ is the K\"{a}hler metric on complex hyperbolic space, $\C H^2$, with constant holomorphic bisectional curvature $-1$. The solution to backward Ricci flow is given by 
\begin{equation}
g(t)=(R^2-6t)g_{\C H^2},\hspace{.2in}-\infty<t<\frac{R^2}{6}.
\end{equation}
\subsection{B9. $S^4$}
Any initial metric can be written as
\begin{equation}
g_0=R^2g_{S^4},
\end{equation}
for some $R>0$. The solution to backward Ricci flow is given by 
\begin{equation}
g(t)=(R^2+6t)g_{S^4},\hspace{.2in}-\frac{R^2}{6}<t<\infty.
\end{equation}
\subsection{B10. $H^4$}
Any initial metric can be written as
\begin{equation}
g_0=R^2g_{H^4}
\end{equation}
for some $R>0$. The solution to backward Ricci flow is given by 
\begin{equation}
g(t)=(R^2-4t)g_{H^4},\hspace{.2in}-\infty<t<\frac{R^2}{6}.
\end{equation}

\section{Conclusions}

Our conclusion is that the end behavior of locally homogeneous manifolds under backward Ricci flow is fairly consistent. In general however, the Bianchi classes have very different, and more interesting, behavior than the non-Bianchi classes. Recall that all of our solutions are Riemannian metrics of the form
\begin{equation*}
g(t)=A(t)\theta_1^2+B(t)\theta_2^2+C(t)\theta_3^2+D(t)\theta_4^2,
\end{equation*}
where we call $A,B,C$ and $D$ the metric coefficients. Each class of manifolds studied in this paper exhibits one or more of the following types of end behavior:\\

Expanding-1: One metric coefficient approaches a constant while the other three grow linearly for all time. The volume-normalized metric approaches a sub-Riemannian geometry, the hyperplane $\R^3$.\\

Expanding-2: Two metric coefficients approach constants while the other two grow linearly for all time. The volume-normalized metric approaches the plane $\R^2$.\\

Expanding-3: All metric coefficients grow linearly for all time. The volume-normalized metric approaches an Einstein metric with positive Ricci curvature.\\

Line-1: Three metric coefficients approach $0$ on the order of $(T_0-t)^{1/3}$ while the other approaches infinity on the order of $(T_0-t)^{-1/3}$. The volume-normalized metric approaches the line $\R$. This is a special case where nearby initial conditions exhibit Pancake-1, Pancake-2 or Pancake-4 end behavior.\\

Line-2: One metric coefficient approaches a constant while the others approach $0$ linearly. The volume-normalized metric approaches the line $\R$.\\

Pancake-1: One metric coefficient approaches $0$ linearly or on the order of $(T_0-t)^{1/3} $, while the other three approach constants. The volume-normalized metric approaches the hyperplane $\R^3$.\\

Pancake-2: One metric coefficient approaches a constant, two approach $0$ on the order of $(T_0-t)^{1/3}$, and the final metric coefficient approaches infinity on the order of $(T_0-t)^{-1/3}$. The volume-normalized metric approaches $\R^2$.\\

Pancake-3: One metric coefficient approaches a constant, one approaches infinity on the order of $(T_0-t)^{-1/3}$, one approaches $0$ on the order of $(T_0-t)^{1/3}$, and the last approaches $0$ on the order of $(T_0-t)^{2/3}$. This is a special case where nearby initial conditions exhibit Pancake-3 end behavior. The volume-normalized metric approaches the plane $\R^2$.\\

Pancake-4: Two metric coefficients approach $0$ linearly while the other two approach constants. The volume-normalized metric approaches the plane $\R^2$.\\

Point: All metric coefficients approach $0$ linearly. Any compact quotients will collapse to a point at $T_0$. Among the classes of manifolds we have considered, this happens only for the Einstein manifolds $\C H^2$ and $H^4$. Thus the volume-normalized metrics are constant metrics.\\

Trivial: Each metric coefficient approaches a constant. Thus, the metric approaches an Einstein metric with constant 0. Among the classes considered, this only happens for trivial metrics which are just quotients of Euclidean Space.\\

Tube: One metric coefficient approaches $0$ linearly while two others approach $0$ at the rate of $(T_0-t)^{1/3}$. The final metric coefficient approaches infinity at a rate of $(T_0-t)^{-1/3}$. This is a special case where nearby initial conditions exhibit Pancake-1 or Pancake-3 end behavior. The volume-normalized metric approaches a metric $g=M^2\times\R$, where $M^2$ is the manifold generated by the two $2$-forms whose coefficients shrink at a rate of $(T_0-t)^{1/3}$.

\begin{table}[p]
\centering
\begin{tabular}{|c|c|c|}
\hline
Lie Algebra Class&Lie Group Structure & Possible End Behaviors\\ \hline\hline
A1. $U1[(1,1,1)]$& $(\R^4,\R^4,\{0\})$ & Trivial\\\hline
A2. $U1[1,1,1]$& $(\tilde{Sol}^3\times\R,\tilde{Sol}^3\times\R,e)$ & Pancake-1\\
& $(Sol^4_0,Sol^4_0,e)$ & Pancake-1\\
& $(Sol^4_{m,n},Sol^4_{m,n},e)$ & Pancake-1\\ \hline
A3. $U1[\Z,\bar{\Z},1]$ & $(\R^4,E(2)\times\R^2,e)$& Pancake-1\\
 & & Pancake-2\\ \hline
A4. $U1[2,1],~\mu=0$ & $(Nil^3\times\R, Nil^3\times\R,e)$ & Pancake-2\\ \hline
A5. $U1[2,1],~\mu=1$ & No Compact Geometries & Pancake-2\\ \hline
A6. $U1[3]$ & $(Nil^4,Nil^4,e)$ & Pancake-2\\
 & & Pancake-3\\ \hline
A7. $U3I0$ & $(Sol^4,Sol^4,e)$ & Line-1\\
 & & Pancake-1 \\
 & & Pancake-2\\
 & & Pancake-3\\
 & & Tube\\ \hline
A8. $U3I2$ & No Compact Geometries & Pancake-2\\
 & & Pancake-3\\\hline
A9. $U3S1$ & $(\tilde{SL}(2,\R)\times\R,\tilde{SL}(2,\R)\times\R,e)$ & Pancake-1\\
 & & Pancake-2\\
 & & Tube\\ \hline
A10. $U3S3$ & $(S^3\times\R,SU(2)\times\R,e)$ & Expanding-1\\
 & & Expanding-2\\
 & & Pancake-2\\ \hline\hline
B1. & $(H^3\times\R,H(3)\times\R,SO(3)\times\{0\})$ & Line-2\\ \hline
B2. & $(S^2\times\R^2,SO(3)\times\R^2,SO(2)\times \{0\})$ & Expanding-2\\ \hline
B3. & $(H^2\times\R^2,S0(3)\times\R^2,SO(2)\times \{0\})$ & Pancake-4\\ \hline
B4. & $(S^2\times S^2,SO(3)\times SO(3), SO(2)\times SO(2))$ & Expanading-3\\ \hline
B5. & $(S^2\times H^2,S0(3)\times H(2),SO(2)\times SO(2))$ & Pancake-4\\ \hline
B6. & $(H^2\times H^2,H(2)\times H(2),SO(2)\times SO(2))$ & Pancake-4\\ \hline
B7. & $(\C P^2,SU(3),U(2))$ & Expanding-3\\ \hline
B8. & $(\C H^2,SU(1,2),U(2))$ & Point\\ \hline
B9. & $(S^4,SO(5),SO(4))$ & Expanding-3\\ \hline
B10. & $(H^4,H(4),SO(4))$ & Point\\ \hline
\end{tabular}
\caption{End Behavior of Backward Ricci Flow}
\end{table}

\newpage

In the non-Bianchi cases, it is clear that the shrinking behavior under backward Ricci flow of negatively curved spaces and the expanding behavior of positively curved spaces is exactly opposite of that seen in forward Ricci flow. Similarly, the behaviors in classes A10 and B1 are reversed.\\

In the Bianchi cases, we notice that in forward Ricci flow all solutions exist for all time except for a certain case of A9: $\hat{SL}(2,\R)\times\R$, where the volume-normalized solution collapses to a plane and exhibits pancake-like behavior. In contrast, the backward Ricci flow produces finite-time singularities  in all cases except for special cases of A10: $SU(2)\times \R$, where the solution increases linearly for all time.\\

It is also worth mentioning that Pancake-1 type behavior in backward Ricci flow occurs as a possibility in exactly the same classes of manifolds which exhibit linear growth in one or more metric coefficient in forward Ricci flow.\\

In generic cases, end-behavior of solutions to the differential equations in forward Ricci flow are of the forms $kt^{\pm 1/3},~k$ or $kt$. Conversely, end-behavior of solutions to the differential equations in backward Ricci flow are of the forms $k(T_0-t)^{\pm 1/3},~k$ or $k(T_0-t)$.

\end{document}